# Explicit design optimization of air rudders for maximizing stiffness and fundamental frequency

Yibo Jia [a], Wen Meng [c], Zongliang Du [a,b,*], Chang Liu [a,b], Shanwei Li [a], Conglei Wang [a],

Zhifu Ge [c], Ruiyi Su [c], Xu Guo [a,b,*]

[a]*State Key Laboratory of Structural Analysis, Optimization and CAE Software for Industrial Equipment, Department of Engineering Mechanics, Dalian University of Technology, Dalian 116023, People's Republic of China*

[b]*Ningbo Institute of Dalian University of Technology, Ningbo 315016, People's Republic of China*

[c]*Beijing System Design Institute of Electro-Mechanic Engineering, Beijing 100854, People's Republic of China*

**Abstract**

In aerospace engineering, there is a growing demand for lightweight design through topology optimization. This paper presents a novel design optimization method for stiffened air rudders, commonly used for aircraft attitude control, based on the Moving Morphable Components (MMC) method. The stiffeners within the irregular enclosed design domain are modeled as MMCs and discretized by shell elements, accurately capturing their geometry and evolution during optimization process using explicit parameters. In order to maximize the stiffness and fundamental frequency of the rudder structures, numerical analysis algorithms were developed with shape sensitivity analysis conducted. To comply with the manufacturing requirement, a minimum thickness is prescribed for the stiffeners. Penalty strategies were developed for the thickness and density of stiffeners with thickness smaller than the threshold to meet the thickness requirement and suppress spurious modes. The method's effectiveness was demonstrated through optimization examples of two typical air rudders, illustrating the significance of stiffener's distribution on design objectives. The explicit modeling characteristics allow for directly importing the optimization results into CAD systems, significantly enhancing the engineering applicability.

*Corresponding author. *E-mail address:* zldu@dlut.edu.cn (Z Du), guoxu@dlut.edu.cn (X Guo)

## 1. Introduction

The pursuit of lightweight design in aircraft structures is a constant endeavor. High-speed aircraft and aerospace vehicles impose greater demands for lightweight designs of components and system-level structures. Topology optimization, seeking to identify the optimal structural layout considering various constraints and design requirements, has garnered significant research attention and witnessed rapid development since the pioneering work of Bendsøe and Kikuchi[1]. In engineering fields[2-4], particularly within the aerospace domain[5-13], various topology methods, e.g., the Solid Isotropic Material with Penalization (SIMP) method[14-16], the Level Set Method (LSM)[17,18], and the Evolutionary Structural Optimization (ESO) methods[19,20] have been applied. In-depth insights into the systematic examination of topology optimization methodologies and recent advancements in the field of aircraft and aerospace structural design can be gained by the references[21-25] and the literatures therein.

Air rudders, serving as attitude control components of aircraft, usually withstand high aerodynamic loads, while a well-designed reinforcing stiffener layout is able to improve the structural performance and achieve weight reduction in aircraft structures. Currently, some prospective studies have conducted optimization desi-

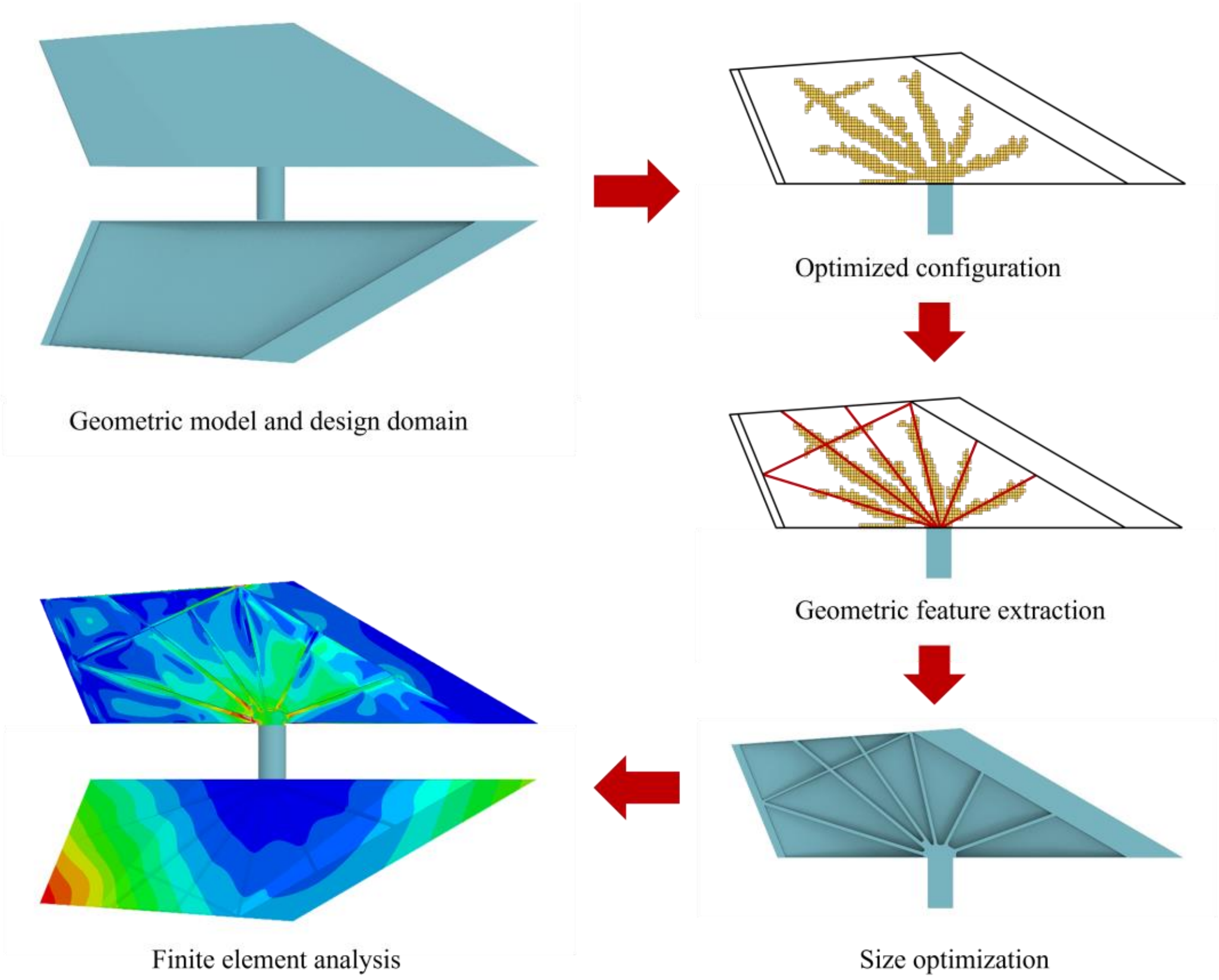

Fig. 1. Traditional optimization design process for rudder structures.

-gn of rudder structures. For instance, Song et al.[10] presented the whole process, from topology optimization design to additive manufacturing of a typical all-movable rudder under thermal-mechanical loads. The aerodynamic shape was meticulously retained in the design to the highest extent possible, and the optimized design underwent rigorous post-processing from an engineering perspective to ensure structural compliance with the constraints of additive manufacturing. Wang et al.[11] presented a multiscale design method using a solid-lattice hybrid structure to enhance the mechanical performance and reduce the structural weight of air rudders, demonstrating a noteworthy enhancement in stiffness compared to the entirely solid and fully lattice configurations. Zhu et al.[12] considered the Y-shaped branches as a unique structural feature. They established a bio-inspired design procedure by integrating the structural layout, sizing parameters, and their simultaneous feature-driven optimization. The benefits were distinctly showcased through its enhanced stiffness and strength compared to the conventional design. Some other studies on the topology optimization of rudder structures is referred to the references[26-28] and literatures therein.

Fig. 1 depicts the traditional optimization design process for air rudders. Implicit optimization method requires discretizing the design domain into solid elements. Therefore, the optimization design of stiffeners for air rudders faces two primary challenges. First, this approach lacks geometric information, making it difficult to ensure that the optimized design results meet the geometric characteristics of stiffeners. The stiffener layout is always manually extracted from the optimized design results, followed by subsequent size optimization to obtain the final design. Secondly, due to the smaller thickness of the stiffener compared to the global size of rudder structures, a dense mesh is necessary to simulate the mechanical behavior accurately. This, in turn, results in increased computational costs and reduced efficiency of the optimization solution process, often serving as a bottleneck that limits the practical implementation of structural topology optimization.

The Moving Morphable Components (MMC) method[29-31] supplies a possible solution approach for improving the issues above. In the MMC method, optimized structures are explicitly described by the geometric parameters of a set of components, which can move, deform, disappear, and overlap during optimization. Because of the explicit topology description characteristics, the MMC method is suitable for optimizing designs involving requirements on geometric features, such as stiffening design of plate and shell structures. Jiang et al. proposed an integrated explicit layout/topology optimization framework[32,33] suitable for designing various thin-walled structures, resulting in reduced computational costs and improved analysis accuracy. They also proposed a method for explicit topology optimization[34] in designing thin-walled structures with complex geometric reinforcing stiffeners based on the technique of computational conformal mapping (CCM). To handle the cases in which shell assumption is not valid, Huo et al. used embedded solid components to construct a general description of thin-walled structures[35], obtaining optimized designs with clear load transmission paths.

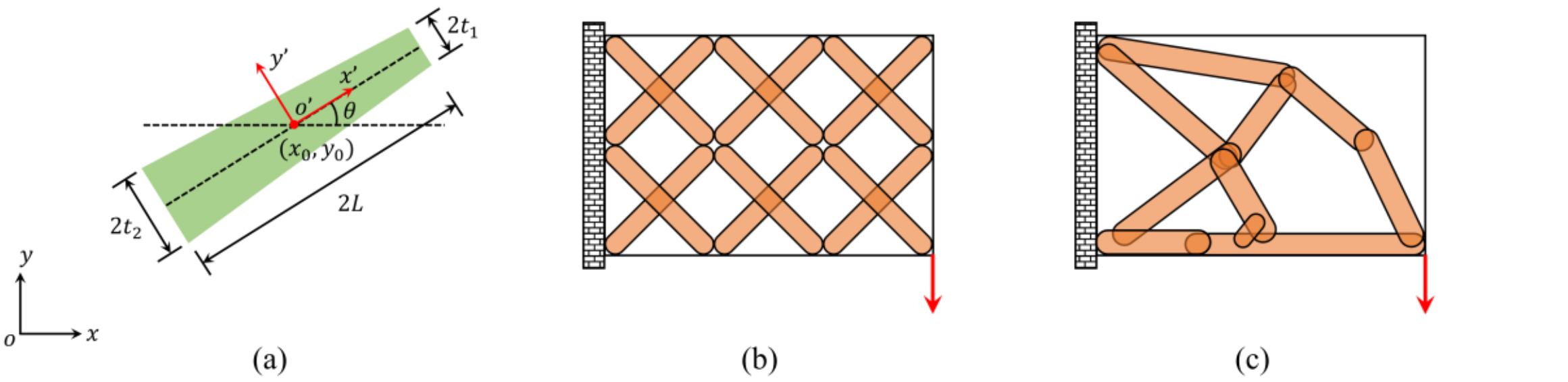

Fig. 2. (a) A Moving Morphable Component with linearly varying thickness; (b) an initial design of the MMC method; (c) the corresponding optimized design.

Nevertheless, those MMC-based design method about thin-walled structures cannot be directly used for the optimal design of air rudder structures mainly due to three aspects: (1) as a closed structures, the height of stiffeners inside the air rudder is varying, while in previous open thin-walled structures, the height of stiffeners is constant; (2) the complexity of air rudder structures poses challenges when using CCM technique under Lagrangian description, in which the mapping error may be non-negligible; (3) the works in references[32-35] focused on the stiffness design problems, while optimization design for dynamic property, e.g., maximizing the fundamental frequency, is not considered.

In this work, the stiffeners are iteratively modelled by the explicit parameters of MMCs, and the Boolean operation is used to obtain the enclosed rudder structure containing irregular stiffeners without using mapping technique. Then shape sensitivities of the structural compliance and fundamental frequency are derived for the stiffeners design of air rudders. To meet the manufacturing requirement about the minimum thickness of stiffeners as well as suppressing spurious modes, penalty strategies are proposed for the thickness and density of stiffeners.

The remaining sections of this paper are organized as follows. Section 2 introduces the explicit description of stiffened rudder structures using the MMC method. Then the mathematical formulation and sensitivity analysis results of static and dynamic optimization of rudder structures are presented in Section 3. In Section 4, two examples of air rudder structures with maximized structural stiffness or fundamental frequency are studied to demonstrate the effectiveness of the proposed method. The conclusions are provided in the last section.

## 2. Explicit description of stiffened air rudders by the Moving Morphable Components (MMC) method

### *2.1. The basic idea of the MMC-based topology optimization*

As shown in Fig. 2, the building blocks of optimized structures are moving, morphable components described by explicit geometric parameters, i.e., the center coordinates, half-length, half-widths, and inclined angle (refer to Fig. 2(a)). By optimizing these geometry parameters, components could move, deform, overlap with each other and disappear, to optimize the structural configuration. Compared to traditional implicit topology optimization methods, the MMC method driven by explicit geometric parameters not only significantly reduces the number of design variables, but also is especially suitable for topology optimization involving geometry requirements, for instance, the additive manufacturing-based topology optimization[36] and topology optimization of thin-walled structures with stiffeners[32-35]. Please refer to the references [29-31] for more details about the MMC method.

### *2.2. Geometry description of structural components of a typical air rudder*

Air rudders are typically triangular or trapezoidal and usually assembled from components such as leading edge, trailing edge, shaft, outer skin and internal stiffeners, as shown in Fig. 3 for a triangular rudder. To satisfy specific practical design requirements, the concept of an enclosed design domain is introduced. The reinforced region is defined as an irregular design domain enclosed by the assembly components, while the structures outside the design domain remain non-designable. Additionally, due to the seamless integration and connection required between the added stiffeners and the upper and lower skins, the height of the stiffeners undergoes continuous variation within the irregular design doma-

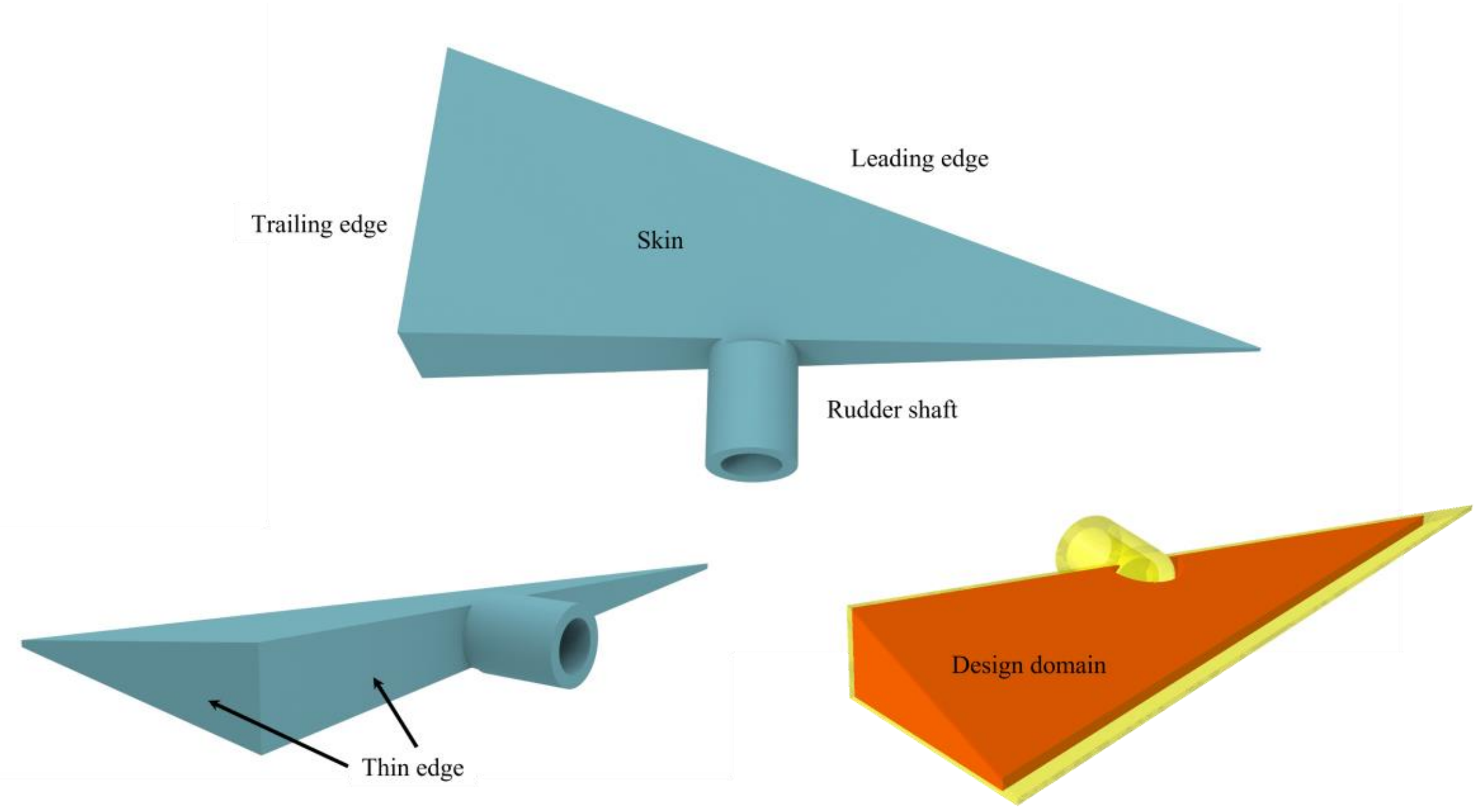

Fig. 3. An illustration of a typical triangular air rudder and its design domain.

-in, posing specific challenges for modeling and optimization tasks.

To enable high-precision numerical analysis while minimizing computational costs, the optimization and analysis models of the stiffened air rudders in this study are formulated using the Lagrangian description. Assuming a reference plane $\mathcal{F}$ located in the $xoy$ plane, as shown in Fig. 4, the $i$-th stiffener within the air rudder can be collectively described by the endpoint coordinates located on the mid-surface ($\boldsymbol{C}^{\alpha}=(C_x^{\alpha},C_y^{\alpha},0)^{\top}$ and $\boldsymbol{C}^{\beta}=(C_x^{\beta},C_y^{\beta},0)^{\top}$) of a linear skeleton $\overline{\boldsymbol{C}^{\alpha}\boldsymbol{C}^{\beta}}$, the stiffener thickness $t^i$, and the height. The height of the $i$-th stiffener on the $\overline{\boldsymbol{C}^{\alpha}\boldsymbol{C}^{\beta}}$ above the reference plane $\mathcal{F}$ $(z \geq 0)$ is denoted as $H_{\mathrm{I}}^i$ (Fig. 4(a)), with its value varying about the positions. Denote the geometric equation of the upper skin plane $\mathcal{M}_{\mathrm{I}}$ above the reference plane $\mathcal{F}$ as:

$$a_{\mathrm{I}}x + b_{\mathrm{I}}y + c_{\mathrm{I}}z + d_{\mathrm{I}} = 0 \tag{1}$$

where $a_{\mathrm{I}}, b_{\mathrm{I}}, c_{\mathrm{I}}, d_{\mathrm{I}}$ represent a set of parameters determined by the geometry of the air rudder. Then the height $H_{\mathrm{I}}^i$ can be represented as follows:

$$H_{\mathrm{I}}^i(\eta) = -\frac{a_{\mathrm{I}}(1-\eta)C_x^{\alpha}}{c_{\mathrm{I}}} - \frac{a_{\mathrm{I}}\eta C_x^{\beta}}{c_{\mathrm{I}}} - \frac{b_{\mathrm{I}}(1-\eta)C_y^{\alpha}}{c_{\mathrm{I}}} - \frac{b_{\mathrm{I}}\eta C_y^{\beta}}{c_{\mathrm{I}}} - \frac{d_{\mathrm{I}}}{c_{\mathrm{I}}} \tag{2}$$

where $\eta \in [0,1]$. By introducing a parameter $\zeta \in [0,1]$ in the height direction, the middle surface $\boldsymbol{\Gamma}_{\mathrm{I}0}^i$ of the $i$-th stiffener in Fig. 4(a) is represented as:

$$\boldsymbol{\Gamma}_{\mathrm{I}0}^i(\eta,\zeta) = \begin{pmatrix} x \\ y \\ z \end{pmatrix} = (1-\eta)\begin{pmatrix} C_x^{\alpha} \\ C_y^{\alpha} \\ 0 \end{pmatrix} + \eta\begin{pmatrix} C_x^{\beta} \\ C_y^{\beta} \\ 0 \end{pmatrix} + \zeta\begin{pmatrix} 0 \\ 0 \\ H_{\mathrm{I}}^i(\eta) \end{pmatrix} \tag{3}$$

The four surfaces comprising the outer boundaries of the $i$-th stiffener depicted in Fig. 4(a) are represented as follows:

$$\boldsymbol{\Gamma}_{\mathrm{I}1}^i(\eta,\zeta) = \boldsymbol{\Gamma}_{\mathrm{I}0}^i(\eta,\zeta) + \frac{t^i}{2}\boldsymbol{n}_{\mathrm{I}}^i \tag{4}$$

$$\boldsymbol{\Gamma}_{\mathrm{I}2}^i(\zeta,\xi) = \boldsymbol{\Gamma}_{\mathrm{I}0}^i(1,\zeta) + \left(\frac{1}{2}-\xi\right)t^i\boldsymbol{n}_{\mathrm{I}}^i \tag{5}$$

$$\boldsymbol{\Gamma}_{\mathrm{I}3}^i(\eta,\zeta) = \boldsymbol{\Gamma}_{\mathrm{I}0}^i(\eta,\zeta) - \frac{t^i}{2}\boldsymbol{n}_{\mathrm{I}}^i \tag{6}$$

$$\boldsymbol{\Gamma}_{\mathrm{I}4}^i(\zeta,\xi) = \boldsymbol{\Gamma}_{\mathrm{I}0}^i(0,\zeta) - \left(\frac{1}{2}-\xi\right)t^i\boldsymbol{n}_{\mathrm{I}}^i \tag{7}$$

Considering the thin geometric characteristics of the stiffener, the minor differences the of surfaces $\boldsymbol{\Gamma}_{\mathrm{I}0}^i, \boldsymbol{\Gamma}_{\mathrm{I}1}^i$, and $\boldsymbol{\Gamma}_{\mathrm{I}3}^i$ are neglected in this work. In Eq. (5) and Eq. (7), $\xi \in [0,1]$ is a positional parameter along the thickness direction, and $\boldsymbol{n}_{\mathrm{I}}^i$ represents the outward normal vector of the boundary $\boldsymbol{\Gamma}_{\mathrm{I}1}^i$.

Similarly, for the portion below the reference plane $\mathcal{F}$ (the $xoy$ plane), the geometric equation of the lower skin plane $\mathcal{M}_{\mathrm{II}}$ is given as:

$$a_{\mathrm{II}}x + b_{\mathrm{II}}y + c_{\mathrm{II}}z + d_{\mathrm{II}} = 0 \tag{8}$$

Other surfaces and geometry parameters of a stiffener below the reference plane are similar as their counterparts of the stiffener above the reference plane, only denoted by a subscript II.

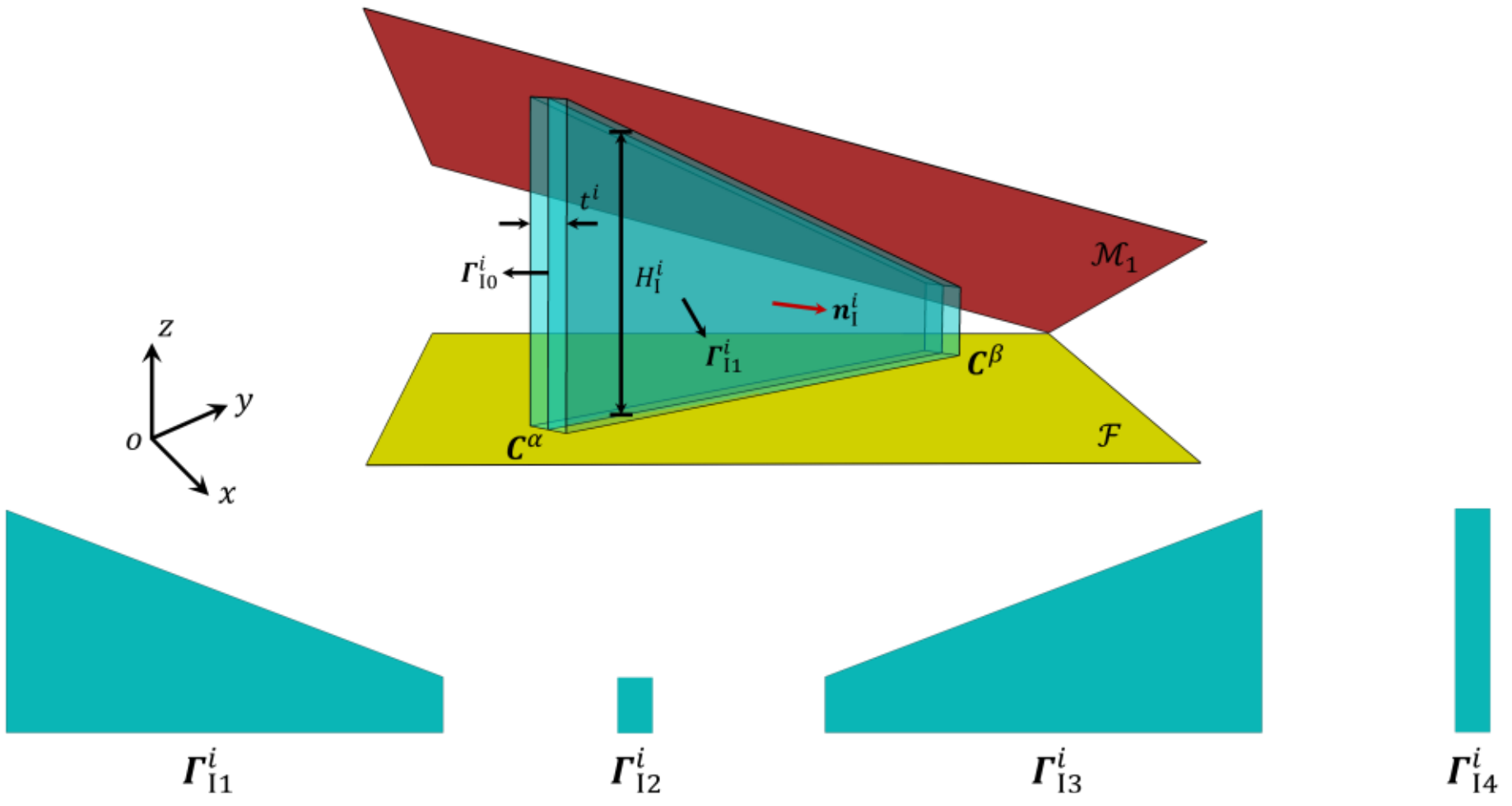

(a) Geometric model and the boundaries of the $i$-th stiffener with $z \geq 0$ (all parameters denoted with subscript I)

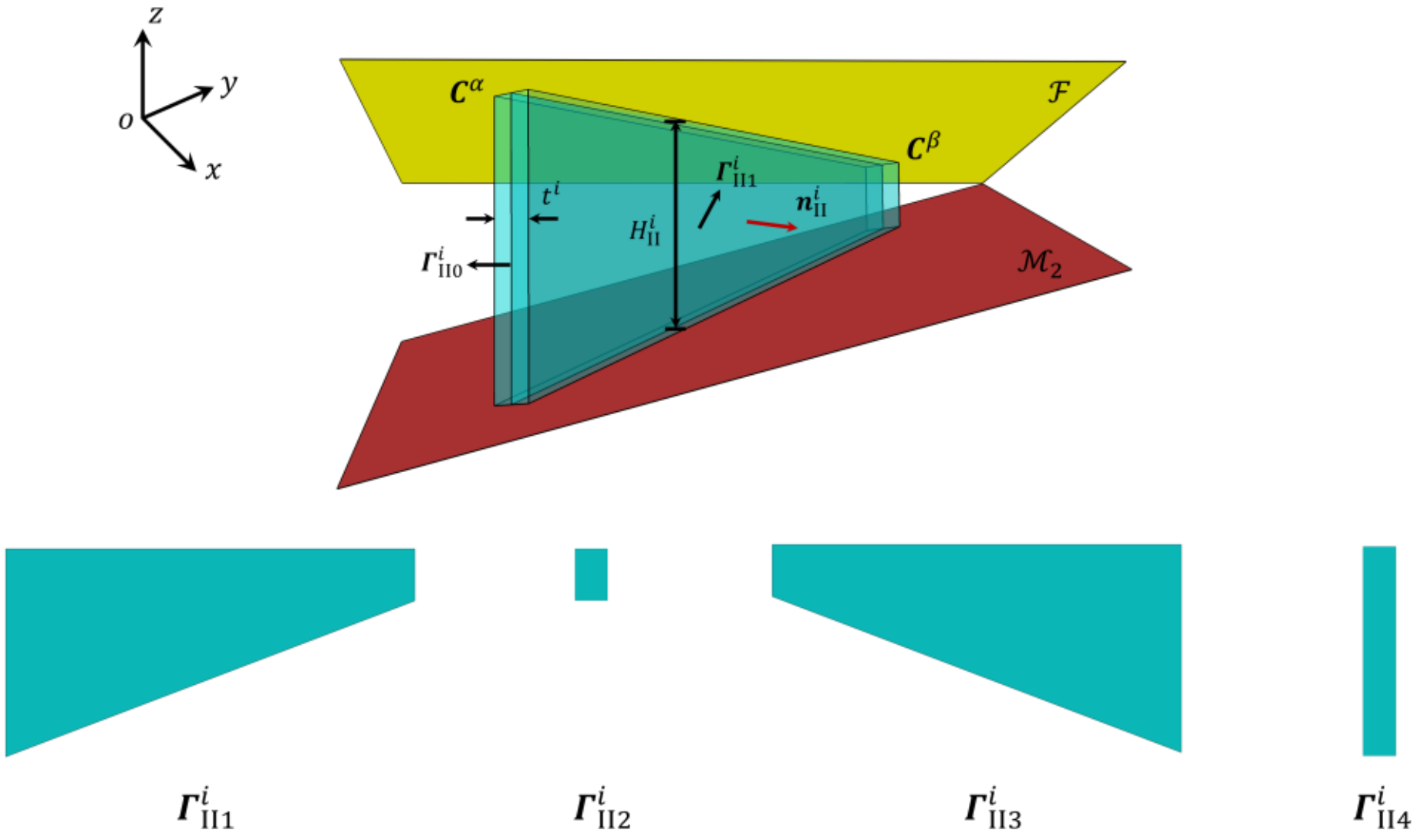

(b) Geometric model and boundaries of the $i$-th stiffener with $z < 0$ (all parameters denoted with subscript II)

Fig. 4. The geometry description of a straight stiffener with a constant thickness.

### *2.3. Modeling and assembly mechanisms*

Based on the Lagrangian description framework, the stiffeners within the enclosed design domain of the air rudder can be explicitly described using geometric parameters. To avoid the intersection and overlap of stiffeners during optimization process, we employ the node-driven adaptive ground structure approach[33] to descript the stiffeners. Compared to the classic ground structure approach, the locations of connection points are also design variables, and the bars are set as stiffeners.

As shown in Fig. 5(a), based on the node-driven adaptive ground structure approach, skeleton of stiffeners is generated within the $xoy$ plane. Then, the skeleton is extruded along the $z$-axis, generating the stiffener layout to be assembled shown in Fig. 5(b). For a typical triangular air rudder, in Fig. 5(c), the stiffeners undergo segmentation using the planes $\mathcal{M}_{\mathrm{I}}, \mathcal{M}_{\mathrm{II}}$, as well as the front plane $\mathcal{M}_{\mathrm{III}}$ for identifying the leading edge. And finally in Fig. 5(d), following the geometric features, the leading edge and the shaft are modeled as solid parts, while the skin and the remaining stiffeners are modeled as shell parts.

It should be noted that, since the exterior of the air rudder are fixed, the above-mentioned modeling and assembly process is valid during the optimization process. Besides, such modeling process ensures a perfect match between the internal stiffeners and skins of air rudders, and is significant for the success of the optimization task.

## 3. Problem formulation and numerical solution aspects

### *3.1. Optimization design for minimizing structural compliance*

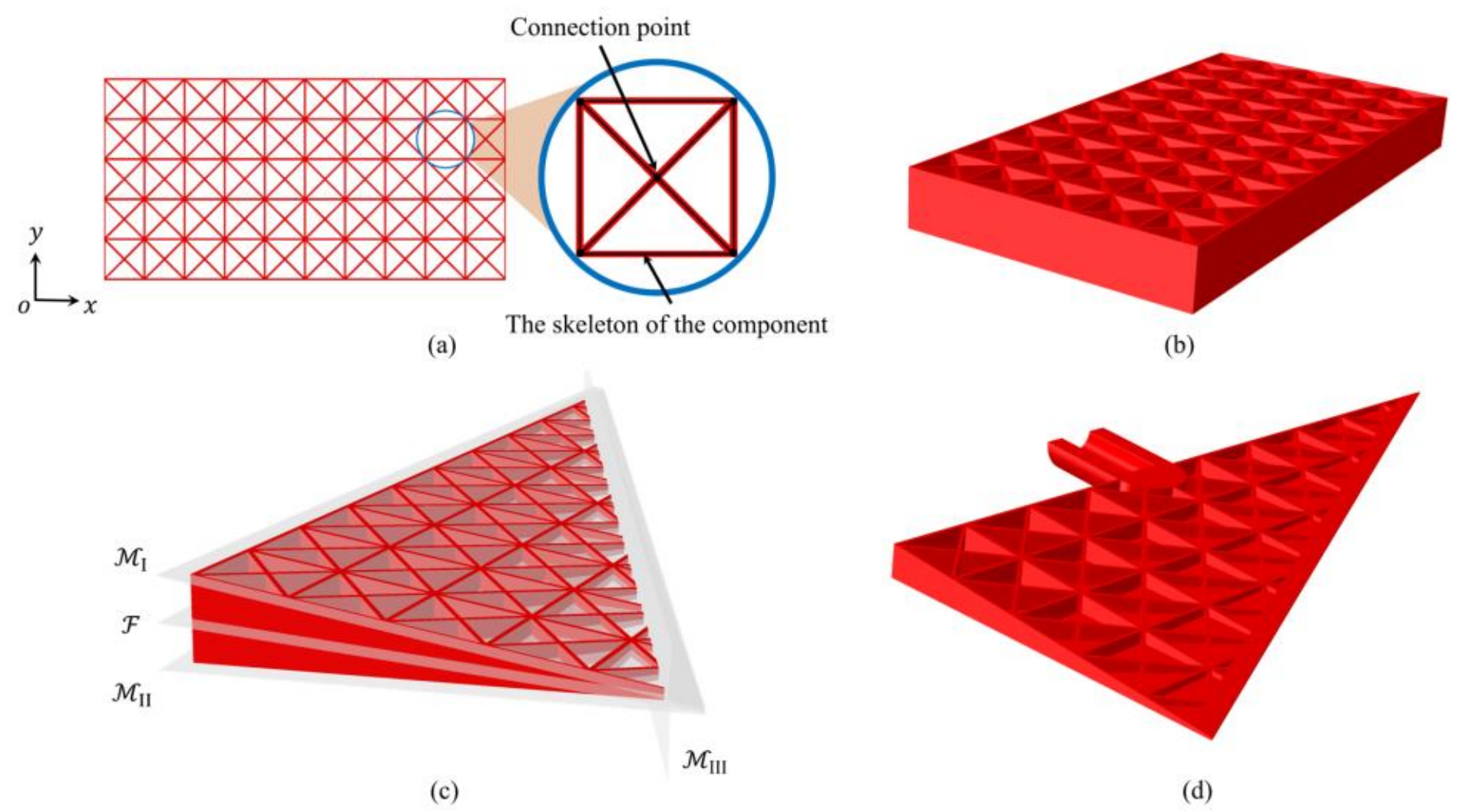

Fig. 5. Schematic diagram of the modeling and assembly mechanism for the air rudder. (a) Skeleton of stiffeners based on the node-driven adaptive ground structure approach; (b) stiffeners extruded from 5(a) to be assembled; (c) cutting the stiffeners by the skins and identifying the leading edge; (d) assembled air rudder (cross-sectional diagram).

Following the node-driven adaptive ground structure approach, it is assumed that all driving nodes are located within the reference surface $\mathcal{F}$. Suppose that the total number of stiffeners within the enclosed design domain is $T$, and the total number of driving nodes is $N$, the design variables $\boldsymbol{I}$ of this problem consist of the thickness parameter vector $\boldsymbol{t} = \left((t^1), \dots, (t^T)\right)^\top$ for stiffeners, and the coordinate vector $\boldsymbol{C} = \left((C_x^1, C_y^1)^\top, \dots, (C_x^N, C_y^N)^\top\right)^\top$ for driving nodes, summarized as $\boldsymbol{I} = (\boldsymbol{t}^\top, \boldsymbol{C}^\top)^\top$. The heights of the $i$-th stiffener $H_{\rm I}^i, H_{\rm II}^i$ vary with the position of driving nodes analytically, hence they are not considered as design variables. The optimization formulation for this problem is presented below:

$$\begin{aligned} &\text{Find} && \boldsymbol{I} \\ &\text{Min.} && C = \boldsymbol{F}^\top \boldsymbol{U} \\ &\text{S.t.} && \mathbf{K}(\boldsymbol{I})\boldsymbol{U}(\boldsymbol{I}) = \boldsymbol{F}, \\ &&& V - \bar{V} \leq 0, \\ &&& \boldsymbol{U} = \overline{\boldsymbol{U}} \quad \text{on } \Gamma_{\boldsymbol{U}}, \\ &&& \boldsymbol{I} \in \mathcal{U}_{\boldsymbol{I}} \end{aligned} \tag{9}$$

where $\boldsymbol{F}$, $\boldsymbol{U}$ and $\mathbf{K}$ respectively represent the nodal force vector, the nodal displacement vector and the global stiffness matrix of the rudder structure. $\bar{V}$ is the upper limit of the volume of stiffeners. $\overline{\boldsymbol{U}}$ represents the prescribed displacement on the Dirichlet boundary $\Gamma_{\boldsymbol{U}}$, and $\mathcal{U}_{\boldsymbol{I}}$ denotes the admissible set of the design variable vector $\boldsymbol{I}$.

### *3.2. Optimization design for maximizing fundamental frequency*

The natural frequencies of an air rudder can be determined by the eigenvalue problem expressed as:

$$(\mathbf{K} - \lambda_i \mathbf{M})\boldsymbol{\varphi}_i = \mathbf{0} \quad i = 1, \dots, n \tag{10}$$

where $\mathbf{K}$ and $\mathbf{M}$ represent the stiffness matrix and mass matrix of the air rudder structure. $\lambda_i = \omega_i^2$ and $\boldsymbol{\varphi}_i$ denote the $i$-th order eigenvalue and the corresponding eigenmode, with $\omega_j$ denoting the $j$-th order natural frequency. $n$ is the total number of degrees of freedom of the structure.

To enhance the dynamic performance of rudder structures, the second optimization problem focuses on maximizing the fundamental frequency of the air rudder as follows:

$$\begin{aligned} &\text{Find} && \boldsymbol{I} \\ &\text{Min.} && G = -\lambda_1 \\ &\text{S.t.} && \left(\mathbf{K}(\boldsymbol{I}) - \lambda_1(\boldsymbol{I})\mathbf{M}(\boldsymbol{I})\right)\boldsymbol{\varphi}_1(\boldsymbol{I}) = \mathbf{0}, \\ &&& V - \bar{V} \leq 0, \\ &&& \boldsymbol{U} = \overline{\boldsymbol{U}} \quad \text{on } \Gamma_{\boldsymbol{U}}, \\ &&& \boldsymbol{I} \in \mathcal{U}_{\boldsymbol{I}} \end{aligned} \tag{11}$$

Unless otherwise specified, the variables in Eq. (11) have the same meaning as those in Eq. (9).

### *3.3. Penalization schemes for stiffener thickness control*

In practical engineering applications, imposing size constraints on structural components holds significant importance in improving design manufacturability. The MMC-based topology optimization framework facilitates direct setting of upper and lower limits for the thickness $t^i$. To ensure that the $i$-th stiffener either is thicker than the manufacturing limit (i.e., $\geq \underline{t}$) or gets eliminated (i.e., $\leq t_\epsilon$) during the optimization pr-ocess, a penalty scheme for thickness is introduced:

$$t_\varepsilon^i = H_\varepsilon^\alpha(t^i - \underline{t})t^i \quad (12)$$

Here, $H_\varepsilon^\alpha(x)$ represents the translated Heaviside function, which can be expressed as:

$$H_\varepsilon^\alpha(x) = \begin{cases} 1 & \text{if } x > \varepsilon, \\ \frac{3(1-\alpha)}{4}\left(\frac{x}{\varepsilon} - \frac{x^3}{3\varepsilon^3}\right) + \frac{1+\alpha}{2} & \text{if } -\varepsilon \le x \le \varepsilon, \\ \alpha & \text{otherwise} \end{cases} \quad (13)$$

where $\varepsilon$ is the regularization parameter and $\alpha$ is a small positive value. In this work, we set $\varepsilon = 0.1$ and $\alpha = 0.001$. Notably, the penalty scheme effectively reduces the number of stiffeners with thickness located at $t_\varepsilon^i \in [t_\epsilon, \underline{t}]$ and guarantees the manufacturability. For more discussions about the thickness penalty scheme, please refer to the reference[33].

In topology optimization problem related to fundamental frequency, a common challenge is suppressing the spurious modes in weak material regions. In particular, for a stiffener, its stiffness is proportional to the cube of thickness while the mass varies linearly. When the thickness decreases, the stiffness decreases faster than the mass, and this would lead to local modes. To address this issue, a density penalization scheme, i.e., the density of $i$-th stiffener material is penalized as well:

$$\rho_\varepsilon^i = \left(H_\varepsilon^\alpha(t^i - \underline{t})\right)^P \rho^i \quad (14)$$

In Eq. (14), $\rho^i$ and $\rho_\varepsilon^i$ represent the original material density and the penalized material density, respectively. To balance the decrease rates of the stiffness and mass of stiffeners, the parameter $P = 2$.

In summary, the thickness penalization scheme ensures the lower bounds of stiffener thickness. The occurrence of spurious modes of thin stiffeners is suppressed by further adoption of the penalization scheme of material density.

### *3.4. Sensitivity analysis*

#### 3.4.1. Sensitivity analysis of structural compliance

The current design method is actually a boundary evolution method based on the Lagrangian description, and it requires the use of a shape sensitivity analysis[37,38] to obtain the derivatives of the objective and constraint functions with respect to the design variables. The weak form of the static equilibrium equation is expressed as:

$$\int_\Omega \mathbb{E}:\boldsymbol{\varepsilon}(\boldsymbol{u}):\boldsymbol{\varepsilon}(\boldsymbol{v})\mathrm{dV} = \int_\Omega \boldsymbol{f}\cdot\boldsymbol{v}\mathrm{dV} + \int_{\Gamma_t} \bar{\boldsymbol{t}}\cdot\boldsymbol{v}\,\mathrm{d}\Gamma, \quad \forall\boldsymbol{v} \in \mathcal{U}_{\mathrm{ad}}^0 \quad (15)$$

where $\Omega$ and $\partial\Omega$ represent the bounded domain and its boundary, respectively. $\mathbb{E}$ is the elasticity tensor of the material, $\boldsymbol{\varepsilon}$ represents the strain tensor, while $\boldsymbol{u}, \boldsymbol{v}$ represent the actual displacement and virtual displacement, respectively. The symbols $\boldsymbol{f}, \bar{\boldsymbol{t}}$ represent the body force density and surface force density acting on the structure. The shape sensitivity of the objective function $C$ in Eq. (9) can be expressed as follows:

$$\delta C = \int_{\partial\Omega} \left(-W(\boldsymbol{u},\boldsymbol{u})\right)V_n \mathrm{d}\Gamma = \sum_i^T \int_{\partial\Omega^i} \left(-W(\boldsymbol{u},\boldsymbol{u})\right)V_n^i \mathrm{d}\Gamma \triangleq \sum_i^T \delta C^i \quad (16)$$

Here, $\partial\Omega^i, i = 1, \dots, T$ represents the boundary of the $i$-th component. $W(\boldsymbol{u},\boldsymbol{u})$ represents the strain energy density, while $V_n^i$ represents the normal velocity field associated with the evolution of the boundary of the $i$-th component. Taking into account the thin geometry characteristics depicted in Fig. 4, the area of the minor faces, denoted as $\boldsymbol{\Gamma}_{\mathrm{I}2}^i, \boldsymbol{\Gamma}_{\mathrm{I}4}^i, \boldsymbol{\Gamma}_{\mathrm{II}2}^i, \boldsymbol{\Gamma}_{\mathrm{II}4}^i$, is smaller compared to $\boldsymbol{\Gamma}_{\mathrm{I}1}^i, \boldsymbol{\Gamma}_{\mathrm{I}3}^i, \boldsymbol{\Gamma}_{\mathrm{II}1}^i, \boldsymbol{\Gamma}_{\mathrm{II}3}^i$. Besides, except for the open end of stiffeners, those minor faces are inside the stiffeners. Therefore, the contribution of those minor faces to the shape sensitivity is neglected, and the shape sensitivity of the $i$-th component can be calculated as follows:

$$\delta C^i \approx \int_{\boldsymbol{\Gamma}_{\mathrm{I}1}^i \cup \boldsymbol{\Gamma}_{\mathrm{I}3}^i \cup \boldsymbol{\Gamma}_{\mathrm{II}1}^i \cup \boldsymbol{\Gamma}_{\mathrm{II}3}^i} \left(-W(\boldsymbol{u},\boldsymbol{u})\right)V_n^i d\Gamma \quad (17)$$

Taking the face $\boldsymbol{\Gamma}_{\mathrm{I}1}^i$ for example, its normal velocity $V_n^i\big|_{\boldsymbol{\Gamma}_{\mathrm{I}1}^i}$ is given as:

$$V_n^i\big|_{\boldsymbol{\Gamma}_{\mathrm{I}1}^i} = \delta\boldsymbol{\Gamma}_{\mathrm{I}1}^i \cdot \boldsymbol{n}_{\mathrm{I}}^i \quad (18)$$

Furthermore, by introducing the variational term for the middle surface $\delta\boldsymbol{\Gamma}_{\mathrm{I}0}^i$ and utilizing Eq. (4) and Eq. (12), we can simplify Eq. (18) as:

$$V_n^i\big|_{\boldsymbol{\Gamma}_{\mathrm{I}1}^i} = \delta\boldsymbol{\Gamma}_{\mathrm{I}0}^i \cdot \boldsymbol{n}_{\mathrm{I}}^i + \frac{\delta t_\varepsilon^i}{2} \quad (19)$$

Since the normal vector $\boldsymbol{n}_{\mathrm{I}}^i$ is perpendicular to the direction vector of the $i$-th stiffener, it can be easily determined that:

$$\boldsymbol{n}_{\mathrm{I}}^i = \left(n_{\mathrm{I}x}^i, n_{\mathrm{I}y}^i, 0\right)^{\mathrm{T}} \quad (20)$$

Where

$$n_{\mathrm{I}x}^i = \frac{C_2^\beta - C_2^\alpha}{\sqrt{\left(C_1^\beta - C_1^\alpha\right)^2 + \left(C_2^\beta - C_2^\alpha\right)^2}} \quad (21)$$

$$n_{\mathrm{I}y}^i = \frac{C_1^\alpha - C_1^\beta}{\sqrt{\left(C_1^\beta - C_1^\alpha\right)^2 + \left(C_2^\beta - C_2^\alpha\right)^2}} \quad (22)$$

So the term $\delta\boldsymbol{\Gamma}_{\mathrm{I}0}^i \cdot \boldsymbol{n}_{\mathrm{I}}^i$ can be expressed as:

$$\delta\boldsymbol{\Gamma}_{\mathrm{I0}}^{i}\cdot\boldsymbol{n}_{\mathrm{I}}^{i}=\delta\left((1-\eta)\boldsymbol{C}^{\alpha}+\eta\boldsymbol{C}^{\beta}+\zeta\left(0,0,H_{\mathrm{I}}^{i}(\eta)\right)^{\mathrm{T}}\right)\cdot\left(n_{\mathrm{I}x}^{i},n_{\mathrm{I}y}^{i},0\right)^{\mathrm{T}}=\left((1-\eta)\delta\boldsymbol{C}^{\alpha}+\eta\delta\boldsymbol{C}^{\beta}+\zeta\left(0,0,\delta H_{\mathrm{I}}^{i}(\eta)\right)^{\mathrm{T}}\right)\cdot\left(n_{\mathrm{I}x}^{i},n_{\mathrm{I}y}^{i},0\right)^{\mathrm{T}} \tag{23}$$

where $\boldsymbol{C}^{\alpha}=(C_x^{\alpha},C_y^{\alpha},0)^{\mathrm{T}}$ and $\boldsymbol{C}^{\beta}=(C_x^{\beta},C_y^{\beta},0)^{\mathrm{T}}$. According to Eq. (2), we have:

$$\delta H_{\mathrm{I}}^{i}(\eta)=-\left(\frac{a_{\mathrm{I}}(1-\eta)}{c_{\mathrm{I}}}\delta C_x^{\alpha}+\frac{a_{\mathrm{I}}\eta}{c_{\mathrm{I}}}\delta C_x^{\beta}+\frac{b_{\mathrm{I}}(1-\eta)}{c_{\mathrm{I}}}\delta C_y^{\alpha}+\frac{b_{\mathrm{I}}\eta}{c_{\mathrm{I}}}\delta C_y^{\beta}\right) \tag{24}$$

Since the outer normal vector $\boldsymbol{n}_{\mathrm{I}}^{i}=\left(n_{\mathrm{I}x}^{i},n_{\mathrm{I}y}^{i},0\right)^{\mathrm{T}}$ lies in the plane $xoy$, the variation term $\delta H_{\mathrm{I}}^{i}(\eta)$ in vector $\left(0,0,\delta H_{\mathrm{I}}^{i}(\eta)\right)^{\mathrm{T}}$ does not contribute to $\delta\boldsymbol{\Gamma}_{\mathrm{I0}}^{i}\cdot\boldsymbol{n}_{\mathrm{I}}^{i}$.

Based on the above facts, the final expression for $V_n^i\big|_{\boldsymbol{\Gamma}_{\mathrm{I1}}^{i}}$ can be obtained as:

$$V_n^i\big|_{\boldsymbol{\Gamma}_{\mathrm{I1}}^{i}}=n_{\mathrm{I}x}^{i}(1-\eta)\delta C_x^{\alpha}+n_{\mathrm{I}x}^{i}\eta\delta C_x^{\beta}+n_{\mathrm{I}y}^{i}(1-\eta)\delta C_y^{\alpha}+n_{\mathrm{I}y}^{i}\eta\delta C_y^{\beta}+\frac{\delta t_{\varepsilon}^{i}}{2} \tag{25}$$

According to Eq. (12), $\delta t_{\varepsilon}^{i}$ can be expanded as:

$$\delta t_{\varepsilon}^{i}=\delta H_{\varepsilon}^{\alpha}\left(t^{i}-\underline{t}\right)t^{i}+H_{\varepsilon}^{\alpha}\left(t^{i}-\underline{t}\right)\delta t^{i} \tag{26}$$

Similarly, considering $\boldsymbol{n}_{\mathrm{II}}^{i}=\boldsymbol{n}_{\mathrm{I}}^{i}=\left(n_{\mathrm{I}x}^{i},n_{\mathrm{I}y}^{i},0\right)^{\mathrm{T}}$, the velocity field terms $V_n^i\big|_{\boldsymbol{\Gamma}_{\mathrm{I3}}^{i}}, V_n^i\big|_{\boldsymbol{\Gamma}_{\mathrm{II1}}^{i}}, V_n^i\big|_{\boldsymbol{\Gamma}_{\mathrm{II3}}^{i}}$ can be represented respectively as:

$$V_n^i\big|_{\boldsymbol{\Gamma}_{\mathrm{I3}}^{i}}=-n_{\mathrm{I}x}^{i}(1-\eta)\delta C_x^{\alpha}-n_{\mathrm{I}x}^{i}\eta\delta C_x^{\beta}-n_{\mathrm{I}y}^{i}(1-\eta)\delta C_y^{\alpha}-n_{\mathrm{I}y}^{i}\eta\delta C_y^{\beta}+\frac{\delta t_{\varepsilon}^{i}}{2} \tag{27}$$

$$V_n^i\big|_{\boldsymbol{\Gamma}_{\mathrm{II1}}^{i}}=n_{\mathrm{I}x}^{i}(1-\eta)\delta C_x^{\alpha}+n_{\mathrm{I}x}^{i}\eta\delta C_x^{\beta}+n_{\mathrm{I}y}^{i}(1-\eta)\delta C_y^{\alpha}+n_{\mathrm{I}y}^{i}\eta\delta C_y^{\beta}+\frac{\delta t_{\varepsilon}^{i}}{2} \tag{28}$$

$$V_n^i\big|_{\boldsymbol{\Gamma}_{\mathrm{II3}}^{i}}=-n_{\mathrm{I}x}^{i}(1-\eta)\delta C_x^{\alpha}-n_{\mathrm{I}x}^{i}\eta\delta C_x^{\beta}-n_{\mathrm{I}y}^{i}(1-\eta)\delta C_y^{\alpha}-n_{\mathrm{I}y}^{i}\eta\delta C_y^{\beta}+\frac{\delta t_{\varepsilon}^{i}}{2} \tag{29}$$

Combining Eqs. (16), (17) and (24)-(29), the sensitivity of the air rudder's structural compliance to the design variables is calculated as follows:

$$\frac{\partial C}{\partial C_x^{\alpha}}=\sum_{i=1}^{T^{\alpha}}\left(\int_{\boldsymbol{\Gamma}_{\mathrm{I1}}^{i}}\left(-W(\boldsymbol{u},\boldsymbol{u})\right)(1-\eta)n_{\mathrm{I}x}^{i}d\Gamma+\int_{\boldsymbol{\Gamma}_{\mathrm{I3}}^{i}}W(\boldsymbol{u},\boldsymbol{u})(1-\eta)n_{\mathrm{I}x}^{i}d\Gamma+\int_{\boldsymbol{\Gamma}_{\mathrm{II1}}^{i}}\left(-W(\boldsymbol{u},\boldsymbol{u})\right)(1-\eta)n_{\mathrm{I}x}^{i}d\Gamma+\int_{\boldsymbol{\Gamma}_{\mathrm{II3}}^{i}}W(\boldsymbol{u},\boldsymbol{u})(1-\eta)n_{\mathrm{I}x}^{i}d\Gamma\right) \tag{30}$$

$$\frac{\partial C}{\partial C_x^{\beta}}=\sum_{i=1}^{T^{\beta}}\left(-\int_{\boldsymbol{\Gamma}_{\mathrm{I1}}^{i}}W(\boldsymbol{u},\boldsymbol{u})\eta n_{\mathrm{I}x}^{i}d\Gamma+\int_{\boldsymbol{\Gamma}_{\mathrm{I3}}^{i}}W(\boldsymbol{u},\boldsymbol{u})\eta n_{\mathrm{I}x}^{i}d\Gamma-\int_{\boldsymbol{\Gamma}_{\mathrm{II1}}^{i}}W(\boldsymbol{u},\boldsymbol{u})\eta n_{\mathrm{I}x}^{i}d\Gamma+\int_{\boldsymbol{\Gamma}_{\mathrm{II3}}^{i}}W(\boldsymbol{u},\boldsymbol{u})\eta n_{\mathrm{I}x}^{i}d\Gamma\right) \tag{31}$$

$$\frac{\partial C}{\partial C_y^{\alpha}}=\sum_{i=1}^{T^{\alpha}}\left(\int_{\boldsymbol{\Gamma}_{\mathrm{I1}}^{i}}\left(-W(\boldsymbol{u},\boldsymbol{u})\right)(1-\eta)n_{\mathrm{I}y}^{i}d\Gamma+\int_{\boldsymbol{\Gamma}_{\mathrm{I3}}^{i}}W(\boldsymbol{u},\boldsymbol{u})(1-\eta)n_{\mathrm{I}y}^{i}d\Gamma+\int_{\boldsymbol{\Gamma}_{\mathrm{II1}}^{i}}\left(-W(\boldsymbol{u},\boldsymbol{u})\right)(1-\eta)n_{\mathrm{I}y}^{i}d\Gamma+\int_{\boldsymbol{\Gamma}_{\mathrm{II3}}^{i}}W(\boldsymbol{u},\boldsymbol{u})(1-\eta)n_{\mathrm{I}y}^{i}d\Gamma\right) \tag{32}$$

$$\frac{\partial C}{\partial C_y^{\beta}}=\sum_{i=1}^{T^{\beta}}\left(-\int_{\boldsymbol{\Gamma}_{\mathrm{I1}}^{i}}W(\boldsymbol{u},\boldsymbol{u})\eta n_{\mathrm{I}y}^{i}d\Gamma+\int_{\boldsymbol{\Gamma}_{\mathrm{I3}}^{i}}W(\boldsymbol{u},\boldsymbol{u})\eta n_{\mathrm{I}y}^{i}d\Gamma-\int_{\boldsymbol{\Gamma}_{\mathrm{II1}}^{i}}W(\boldsymbol{u},\boldsymbol{u})\eta n_{\mathrm{I}y}^{i}d\Gamma+\int_{\boldsymbol{\Gamma}_{\mathrm{II3}}^{i}}W(\boldsymbol{u},\boldsymbol{u})\eta n_{\mathrm{I}y}^{i}d\Gamma\right) \tag{33}$$

$$\frac{\partial C}{\partial t^{i}}=\int_{\boldsymbol{\Gamma}_{\mathrm{I1}}^{i}\cup\boldsymbol{\Gamma}_{\mathrm{I3}}^{i}\cup\boldsymbol{\Gamma}_{\mathrm{II1}}^{i}\cup\boldsymbol{\Gamma}_{\mathrm{II3}}^{i}}-\frac{1}{2}\left(\frac{\partial H_{\varepsilon}^{\alpha}\left(t^{i}-\underline{t}\right)}{\partial t^{i}}t^{i}+H_{\varepsilon}^{\alpha}\left(t^{i}-\underline{t}\right)\right)W(\boldsymbol{u},\boldsymbol{u})d\Gamma \tag{34}$$

where $T^{\alpha}$ and $T^{\beta}$ represent the number of stiffeners associated with the driving points $\boldsymbol{C}^{\alpha}$ and $\boldsymbol{C}^{\beta}$, respectively. From Eq. (13), the exact expression of $\frac{\partial H_{\varepsilon}^{\alpha}(t^{i}-\underline{t})}{\partial t^{i}}$ is trivial and will not be presented.

3.4.2. Sensitivity analysis of structural fundamental frequency

The weak form of governing equation for the fundamental frequency problem is:

$$\int_{\Omega} \mathbb{E} : \boldsymbol{\varepsilon}(\boldsymbol{u}) : \boldsymbol{\varepsilon}(\boldsymbol{v}) \mathrm{dV} = \lambda_1 \int_{\Omega} \rho_{\varepsilon}^{i} \boldsymbol{u} \cdot \boldsymbol{v} \mathrm{dV}, \quad \forall \boldsymbol{v} \in \mathcal{U}_{\mathrm{ad}}^{0} \tag{35}$$

Whereas $\lambda_1$ represents the non-repeated fundamental eigenvalue. For the repeated eigenvalue problem in structural optimization, please refer to reference[39], which is not considered in this study. Since this optimization problem is also self-adjoint[40-42], the shape sensitivity of the objective function $G$ in Eq. (11) can be expressed as:

$$\delta G \approx \sum_{i}^{T} \int_{\Gamma_{\mathrm{I1}}^{i} \cup \Gamma_{\mathrm{I3}}^{i} \cup \Gamma_{\mathrm{II1}}^{i} \cup \Gamma_{\mathrm{II3}}^{i}} \left( \lambda_1 \rho_{\varepsilon}^{i} \boldsymbol{u} \cdot \boldsymbol{u} - W(\boldsymbol{u}, \boldsymbol{u}) \right) V_n^i \mathrm{d}\Gamma \tag{36}$$

By replacing the term $-W(\boldsymbol{u}, \boldsymbol{u})$ in Eqs. (30)-(34) as $\lambda_1 \rho_{\varepsilon}^{i} \boldsymbol{u} \cdot \boldsymbol{u} - W(\boldsymbol{u}, \boldsymbol{u})$, the sensitivity of the fundamental frequency of the air rudder to the design variables can be obtained.

### 3.4.3. Sensitivity analysis of structural volume

The volume $V^i$ of the $i$-th stiffener can be analytically expressed as:

$$V^i = L^i H_*^i t_{\varepsilon}^i \tag{37}$$

where $L^i$ represents the length of the stiffener and $H_*^i$ represents the height of the stiffener at the midpoint of the $\overline{\boldsymbol{C}^{\alpha}\boldsymbol{C}^{\beta}}$. They can be respectively expressed as:

$$L^i = \sqrt{\left(C_1^{\beta} - C_1^{\alpha}\right)^2 + \left(C_2^{\beta} - C_2^{\alpha}\right)^2} \tag{38}$$

$$H_*^i = H_{\mathrm{I}}^i(\eta) - H_{\mathrm{II}}^i(\eta), \eta = 0.5 \tag{39}$$

Therefore, the variation of $V^i$ can be calculated as:

$$\delta V^i = \delta L^i H_*^i t_{\varepsilon}^i + L^i \delta H_*^i t_{\varepsilon}^i + L^i H_*^i \delta t_{\varepsilon}^i \tag{40}$$

The sensitivities of the volume constraint are detailed as:

$$\frac{\partial V}{\partial C_x^{\alpha}} = \sum_{i=1}^{T^{\alpha}} -\frac{\left(C_x^{\beta} - C_x^{\alpha}\right)}{L^i} H_*^i t_{\varepsilon}^i + \frac{1}{2}\left(\frac{a_{\mathrm{II}}}{c_{\mathrm{II}}} - \frac{a_{\mathrm{I}}}{c_{\mathrm{I}}}\right) L^i t_{\varepsilon}^i \tag{41}$$

$$\frac{\partial V}{\partial C_x^{\beta}} = \sum_{i=1}^{T^{\beta}} \frac{\left(C_x^{\beta} - C_x^{\alpha}\right)}{L^i} H_*^i t_{\varepsilon}^i + \frac{1}{2}\left(\frac{a_{\mathrm{II}}}{c_{\mathrm{II}}} - \frac{a_{\mathrm{I}}}{c_{\mathrm{I}}}\right) L^i t_{\varepsilon}^i \tag{42}$$

$$\frac{\partial V}{\partial C_y^{\alpha}} = \sum_{i=1}^{T^{\alpha}} -\frac{\left(C_y^{\beta} - C_y^{\alpha}\right)}{L^i} H_*^i t_{\varepsilon}^i + \frac{1}{2}\left(\frac{b_{\mathrm{II}}}{c_{\mathrm{II}}} - \frac{b_{\mathrm{I}}}{c_{\mathrm{I}}}\right) L^i t_{\varepsilon}^i \tag{43}$$

$$\frac{\partial V}{\partial C_y^{\beta}} = \sum_{i=1}^{T^{\beta}} \frac{\left(C_y^{\beta} - C_y^{\alpha}\right)}{L^i} H_*^i t_{\varepsilon}^i + \frac{1}{2}\left(\frac{b_{\mathrm{II}}}{c_{\mathrm{II}}} - \frac{b_{\mathrm{I}}}{c_{\mathrm{I}}}\right) L^i t_{\varepsilon}^i \tag{44}$$

$$\frac{\partial V}{\partial t^i} = L^i H_*^i \left( \frac{\partial H_{\varepsilon}^{\alpha}\left(t^i - \underline{t}\right)}{\partial t^i} t^i + H_{\varepsilon}^{\alpha}\left(t^i - \underline{t}\right) \right) \tag{45}$$

## 4. Numerical examples

In this section, we aim to demonstrate the effectiveness of the proposed method for optimal design of air rudders, aiming at minimizing structural compliance and maximizing fundamental frequency, respectively. The Method of Moving Asymptotes (MMA) is used as the optimizer[43]. Under the premise of $(V^{(k)} - \bar{V})/\bar{V} \leq 10^{-2}$, the iteration process terminates once the condition $\|C^{(k)} - C^{(k-1)}\| / \|C^{(k)}\| \leq 10^{-2}$ or $\|G^{(k)} - G^{(k-1)}\| / \|G^{(k)}\| \leq 10^{-2}$ is achieved for five consecutive iteration steps. For the manufacturing requirement, the lower bound of stiffeners is set as $\underline{t} = 4.0$mm.

### *4.1. Structural compliance minimization of a trapezoidal air rudder*

The dimensions of the trapezoidal air rudder are illustrated in Fig. 6. The rudder structure is symmetric about the plane $\mathcal{F}$, with the skin thickness of 1.2mm. The design domain of stiffeners is the internal irregular space covered by the leading edge, skin and shaft in Fig. 6., and the material properties listed in Table 1. To approximately simulate the loading conditions during flight, a uniform pressure with amplitude of 0.2MPa is applied to the top skin of the air rudder. Similarly, the amplitude of the pressure on the lower skin is 0.001MPa. Besides, outer cylindrical surface of the shaft is fully fixed. The admissible volume fraction in the enclosed design domain should be no greater than 33%.

Table 1 Material properties of the trapezoidal air rudder.

| | |
|---|---|
| Density $\rho$(Kg/m$^3$) | 4.45E3 |
| Young's modulus $E$(GPa) | 100 |
| Poisson's ratio $\nu$ | 0.3 |

The initial layout of the stiffeners is illustrated in Fig. 7, which comprises a total of 189 stiffeners. After applying the penalty method as described in Eq. (12), the thickness of stiffeners located on the boundaries is $t_{\varepsilon} = 5$mm (the corresponding design variable $t = 5$mm), while the rest have a uniform thickness of $t_{\varepsilon} = 0.002$mm (the corresponding design variable $t = 2$mm). The design variables consist of the stiffener thickness within the interval of [0.001,15]mm, the movement limits of the $x$-coordinate and the $y$-coordinate of the driving nodes are within [−8.775,8.775]mm and [−10.27,10.27]mm at each iteration, respectively.

By solving the formulation Eq. (9), the optimization iteration history curves for the objective function and volume fraction are presented in Fig. 8(a). The structural compliance experiences a rapid decrease within the first 20 iterations, followed by a gradual decline. The volume fraction approaches its upper bound after 107 iterations, and the structural compliance converges to 2452.61mJ with total weight of 4.02 kg at the 185 iteration. Thanks to the thickness penalization scheme, thin stiffeners that fail to meet the thickness constraint have less impact on the structural compliance. After removing the stiffeners with thickness in the interval of $[0.001, 4]$mm, the final optimized result satisfying manufacturing constraint is depicted in Fig. 8(b), of which the final total mass and structural compliance are 4.01 kg and 2491.64 mJ, respectively.

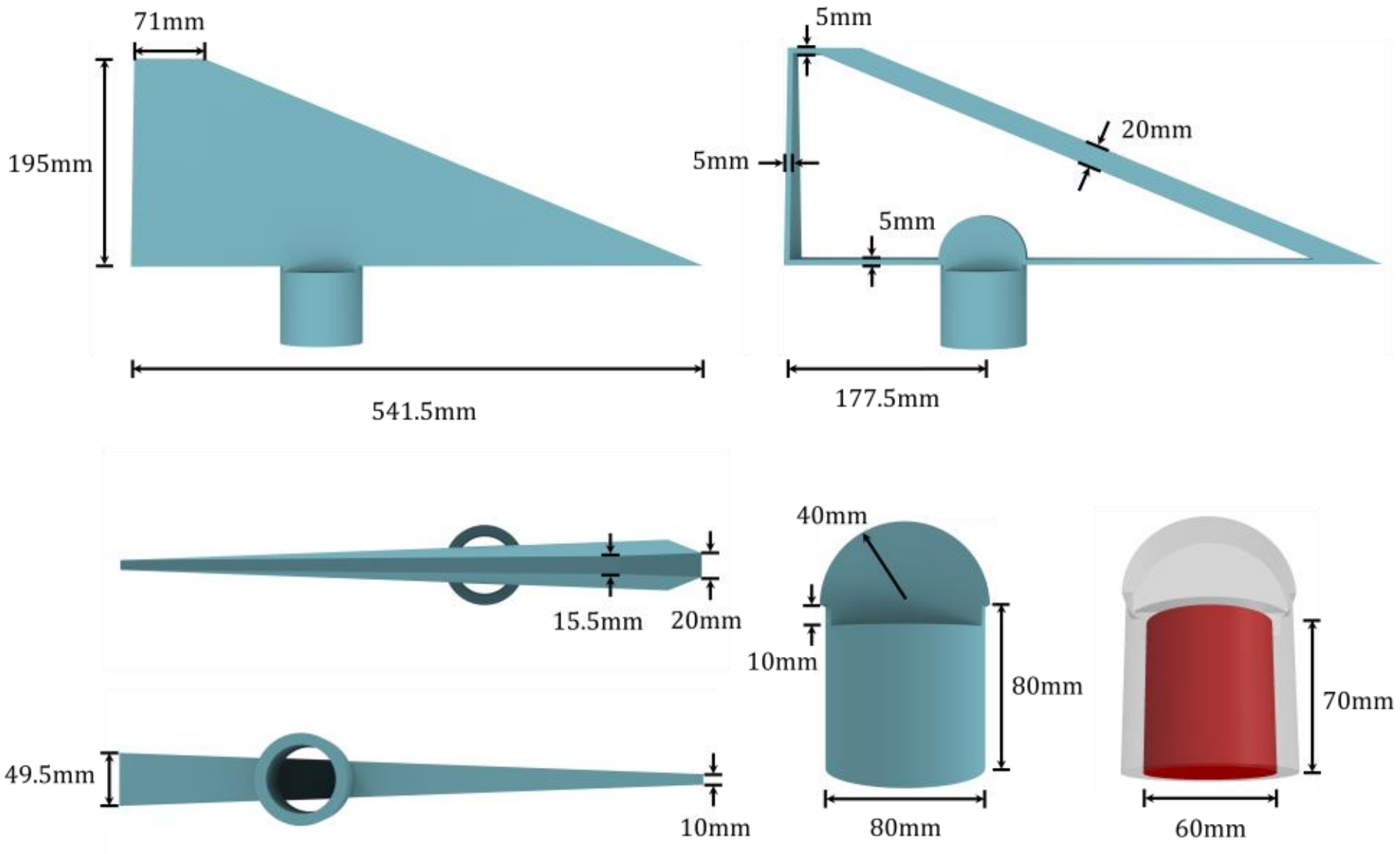

Fig. 6. Dimensions of a trapezoidal air rudder (red cylinder illustrate the void region).

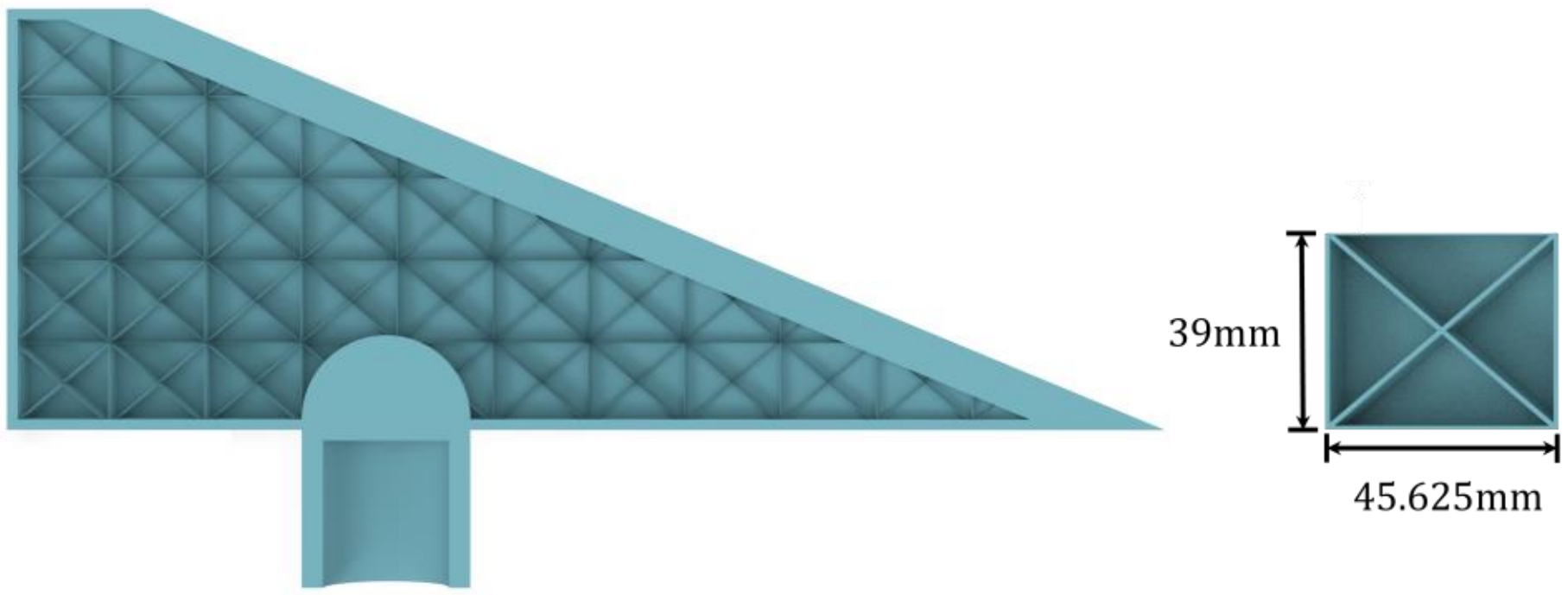

Fig. 7. Initial layout of stiffeners displaying thickness and a primitive cell.

The design domain is divided into 4 regions as Fig. 9(a). One of the advantages of the present design method is that the optimized design has explicit geometry description and it can be modeled in CAD software without tedious postprocessing. In Fig. 9(b), the thicknesses of the stiffeners in Region $R_1$ is illustrated, and there are five stiffeners "growing" from the rudder shaft to supply load transmission path to Regions $R_2, R_3$ and $R_4$. We also noticed that, the thicknesses of the optimized stiffeners are relatively close. The average thickness of the stiffeners in the four regions are $4.29\text{mm}, 4.06\text{mm}, 4.08\text{mm}$ and $4.09\text{mm}$, respectively. This is because uniform pressures are applied on the skin, and uniformly distributed stiffeners are better to suppress the overall deformation.

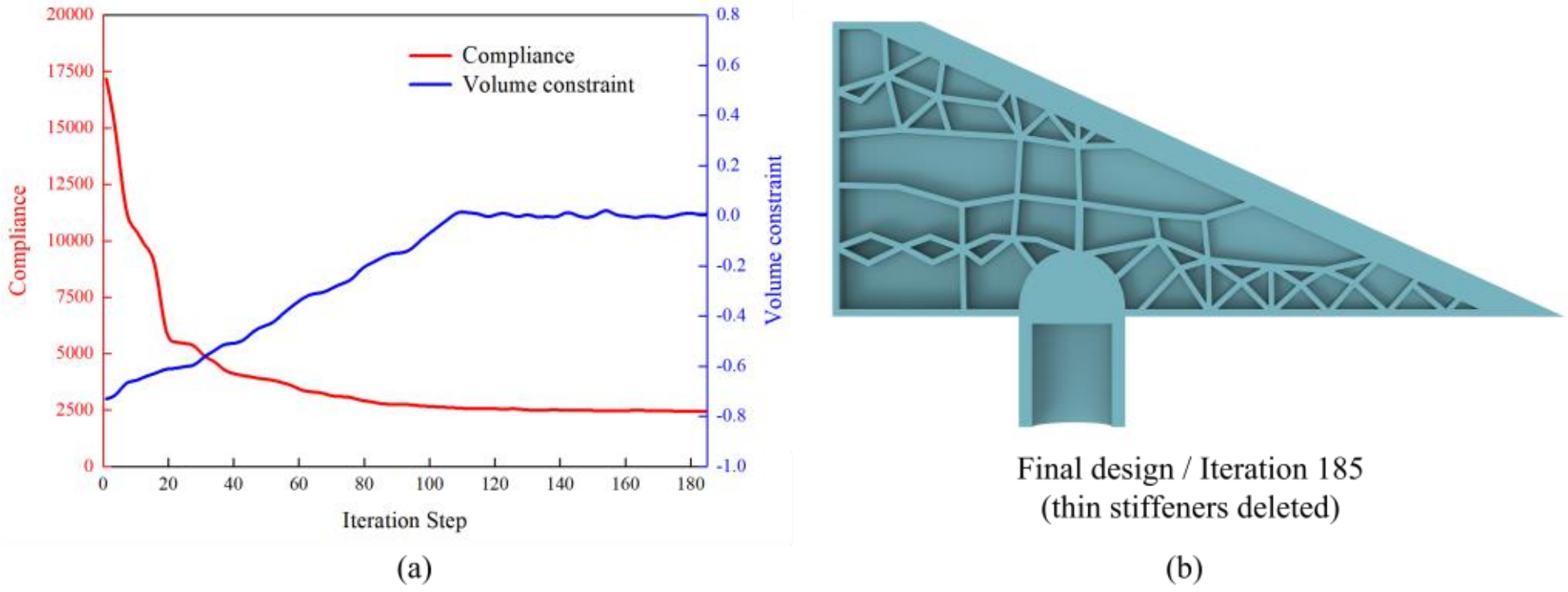

Fig. 8. (a) Iteration history of objective function and constraint values; (b) final design of trapezoidal rudder displaying thickness.

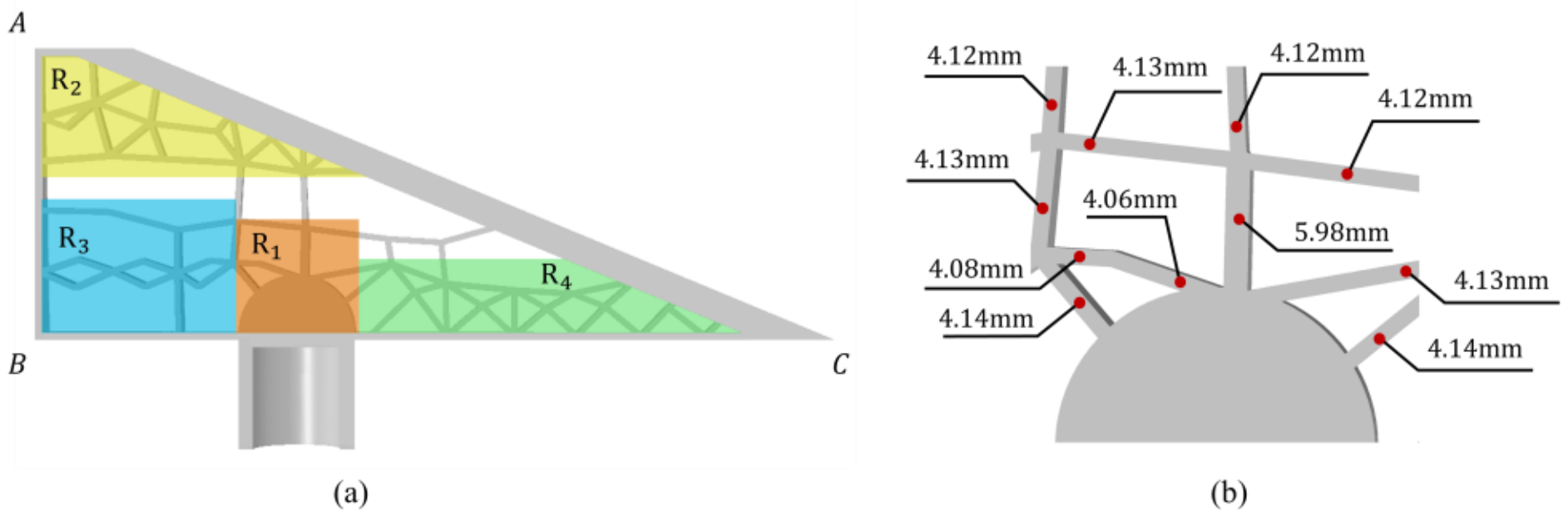

Fig. 9. (a) Distribution of stiffeners concentrated area; (b) thickness of stiffeners in $R_1$.

To better illustrate the optimization process, the intermediate stiffener distribution (displaying thickness), stress field, and displacement field at the 10th, 20th, 80th, and 185th iterations are presented in Table 2. At the 10th iteration, the material concentration is predominantly observed at the shaft and corner point $A$. However, due to the too small thickness of the remaining stiffeners, the deformation is primarily characterized by the large inward deflection of the skin and the downward displacement of corner point $C$. This deformation pattern leads to high stress levels. At the 20th iteration, the stiffeners along the shaft connect to the leading edge, forming a load transformation path, which results in a significant reduction on the maximum displacement amplitude (3.25mm vs 3.02mm) and global stress level decreases compared to the 10th iteration. At the 80th iteration, thicker stiffeners appear in regions $R_2, R_3, R_4$ as well, leading to further improvement on the structural stiffness. In the final configuration at the 185th iteration (without removing thin stiffeners), the maximum displacement amplitude is significantly decreased by 26.46% compared to the design at 10th iteration (2.39mm vs 3.25mm), and the compliance of the whole structure has been decreased by 76.84%. By removing the thin stiffeners in the optimized result, the stress field and displacement field show no significant changes compared to the previous results (maximum displacement amplitude: 2.39mm vs 2.39mm), the stress levels in both the stiffeners and the skin are all below 161.28MPa, validating the effectiveness of the thickness penalization scheme.

### *4.2. Maximizing the structural fundamental frequency of a triangular air rudder*

In the second example, the stiffeners of a single-shaft triangular air rudder symmetric along the plane $\mathcal{F}$ are optimized to maximize the fundamental frequency. The dimensions of the air rudder are illustrated in Fig. 10 with the skin thickness of 1.5mm, and the material properties listed in Table 3. The outer cylindrical surface of the shaft is fully fixed, and the admissible volume fraction in the enclosed design domain should be no greater than 38%. As shown in Fig. 11(a), the initial layout comprises a total of 154 stiffeners. After applying the penalty method as described in Eq. (12), there have a uniform thickness of $t_\varepsilon = 0.0035\text{mm}$ (the corresponding design variable $t = 3.5\text{mm}$) beside boundaries.

Table 3 Material properties of the triangular air rudder.

| | |
|---|---|
| Density $\rho(\text{Kg/m}^3)$ | 4.50E3 |
| Young's modulus $E(\text{GPa})$ | 88 |
| Poisson's ratio $\nu$ | 0.39 |

To validate the effectiveness of the density penalization approach in suppressing spurious modes, Fig. 11(b) and Fig. 11(c) present the vibration modes corresponding to the fundamental frequency for the layout depicted in Fig. 11(a). In the absence of density penalization, the structure exhibits a highly localized spurious mode with a frequency of $f_1 = 13.2\text{Hz}$ as shown in Fig. 11(b). By employing the density penalization scheme using Eq. (14), the fundamental frequency is $f_1 = 191.5\text{Hz}$, and the vibration mode is a global pattern of out-of-plane bending as Fig. 11(c).

The iteration curves of the fundamental frequency and volume constraint are presented in Fig. 12(a). The fundamental frequency of the structure exhibits a sharp increase before the 11th iteration then gradually rises thereafter, and eventually converges at the 171th iteration. Fig. 12(a) also shows the total volume of the thin stiffeners ($t_\varepsilon^i < 4.0\text{mm}$) is relatively larger at the initial stage (maximum 4.8% of the upper limit of volume), and this value decreases to 2.6% finally. Although those thinner stiffeners are modelled during the whole optimization process, the fundamental freq-

Table 2 Configurations and analysis results for some intermediate and final designs.

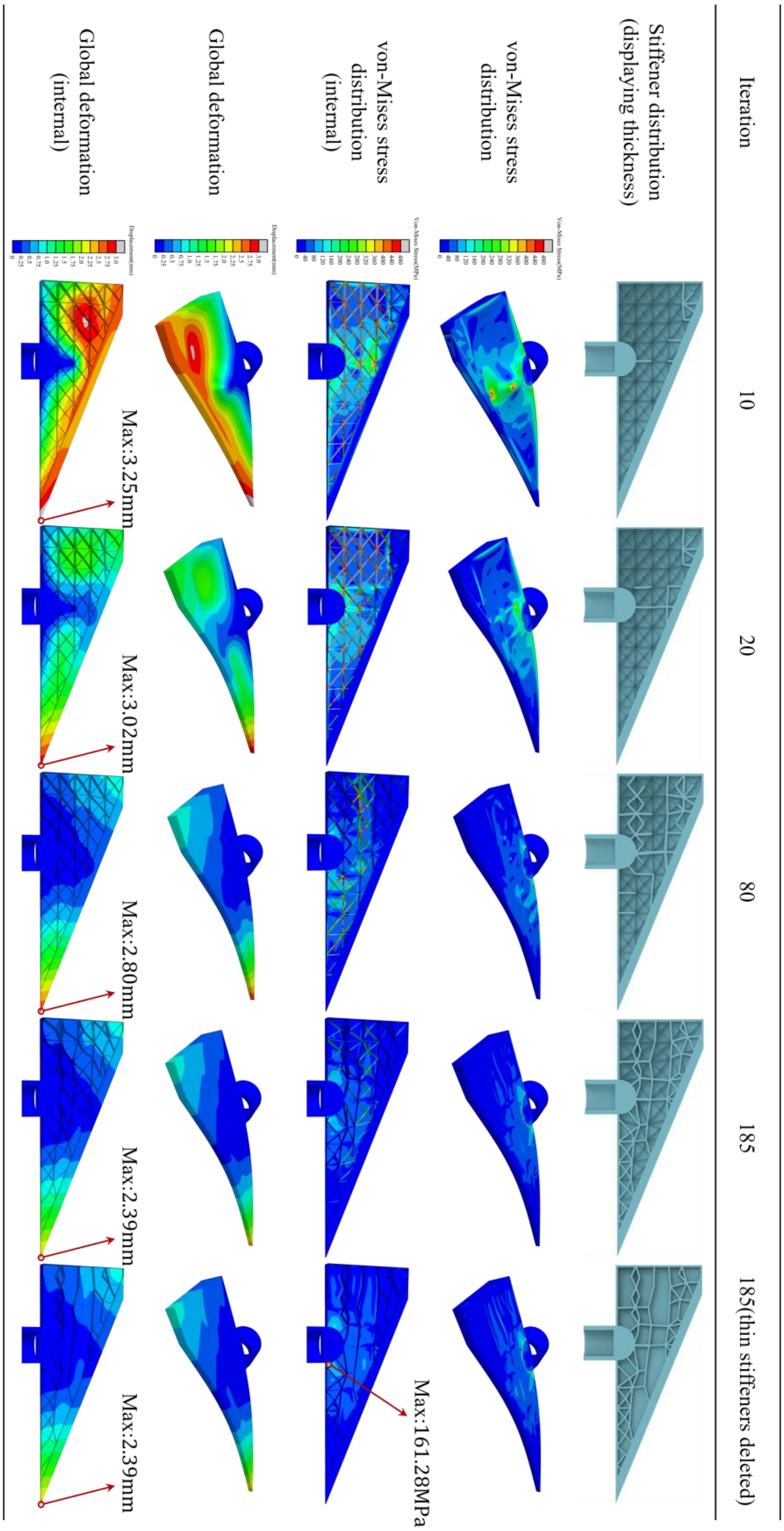

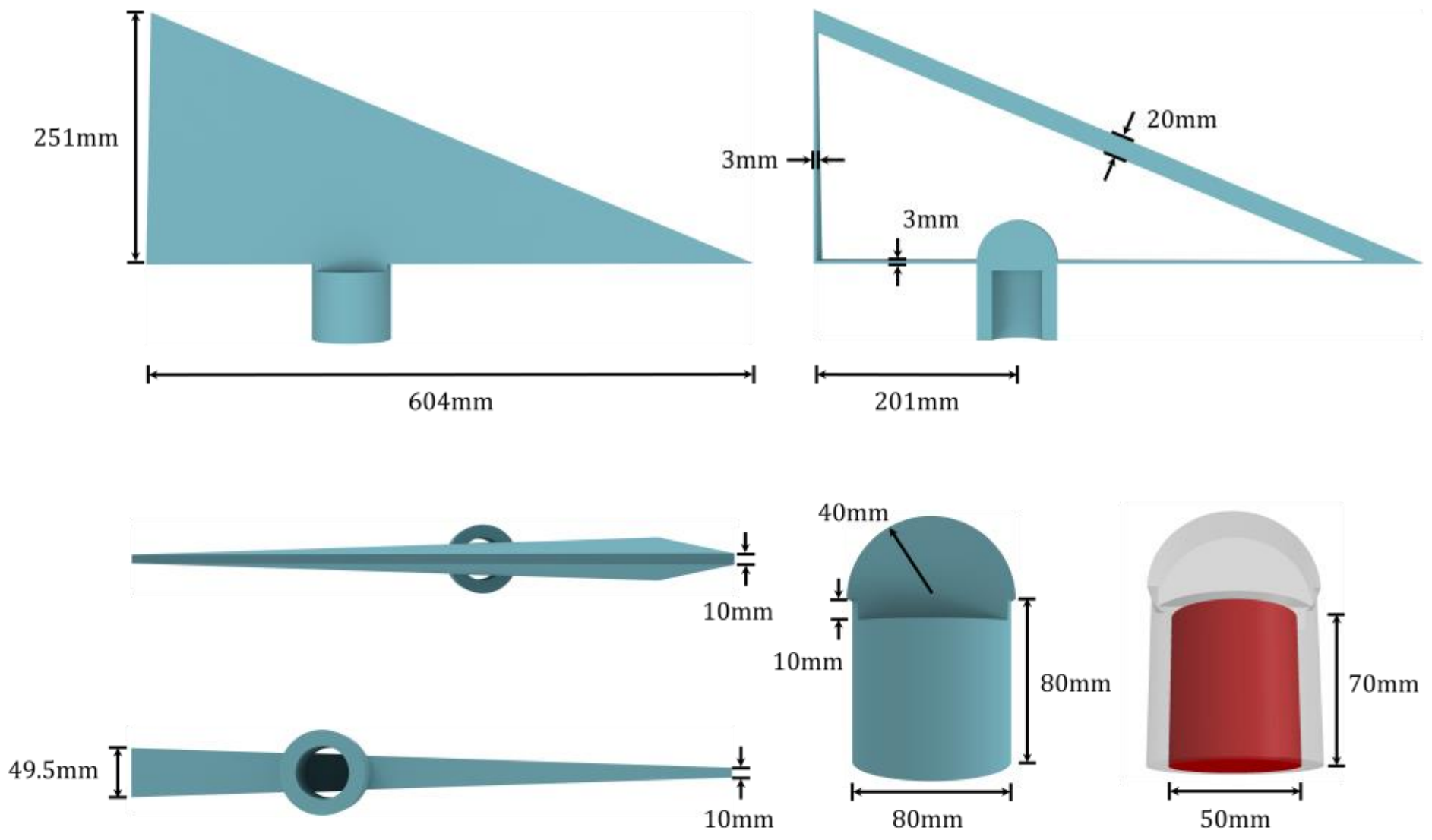

Fig. 10. Dimensions of a triangular air rudder (red cylinder illustrating the void region).

-uency keeps stable and no spurious mode effect is observed. By removing the thin stiffeners with thickness smaller than $4.0\mathrm{mm}$, the final optimized design is depicted in Fig. 12(b), which has a total mass of $5.26\,\mathrm{kg}$ and a fundamental frequency of $f_1 = 249.2\mathrm{Hz}$. The stiffeners are thicker close to the constrained shaft and gradually become thinner gradually outward. This is quite reasonable in mechanics and consistent with the optimized designs in references[10-12].

Furthermore, a typical triangular air rudder was chosen as a reference design to comparatively validate the effectiveness of the optimization method, and the analysis results are listed in Table 4. The original design possesses the same mass as the optimized design, with an out-of-plane bending vibration mode and a fundamental frequency of $217.0\mathrm{Hz}$. The optimized design demonstrates a $14.84\%$ increase in the fundamental frequency compared to the reference, resulting a significant improvement.

To evaluate the optimization process, Table 5 presents the fundamental frequency of intermediate designs at iterations $11, 80$ and $171$, respectively. By modelling and analyzing those designs and their co-

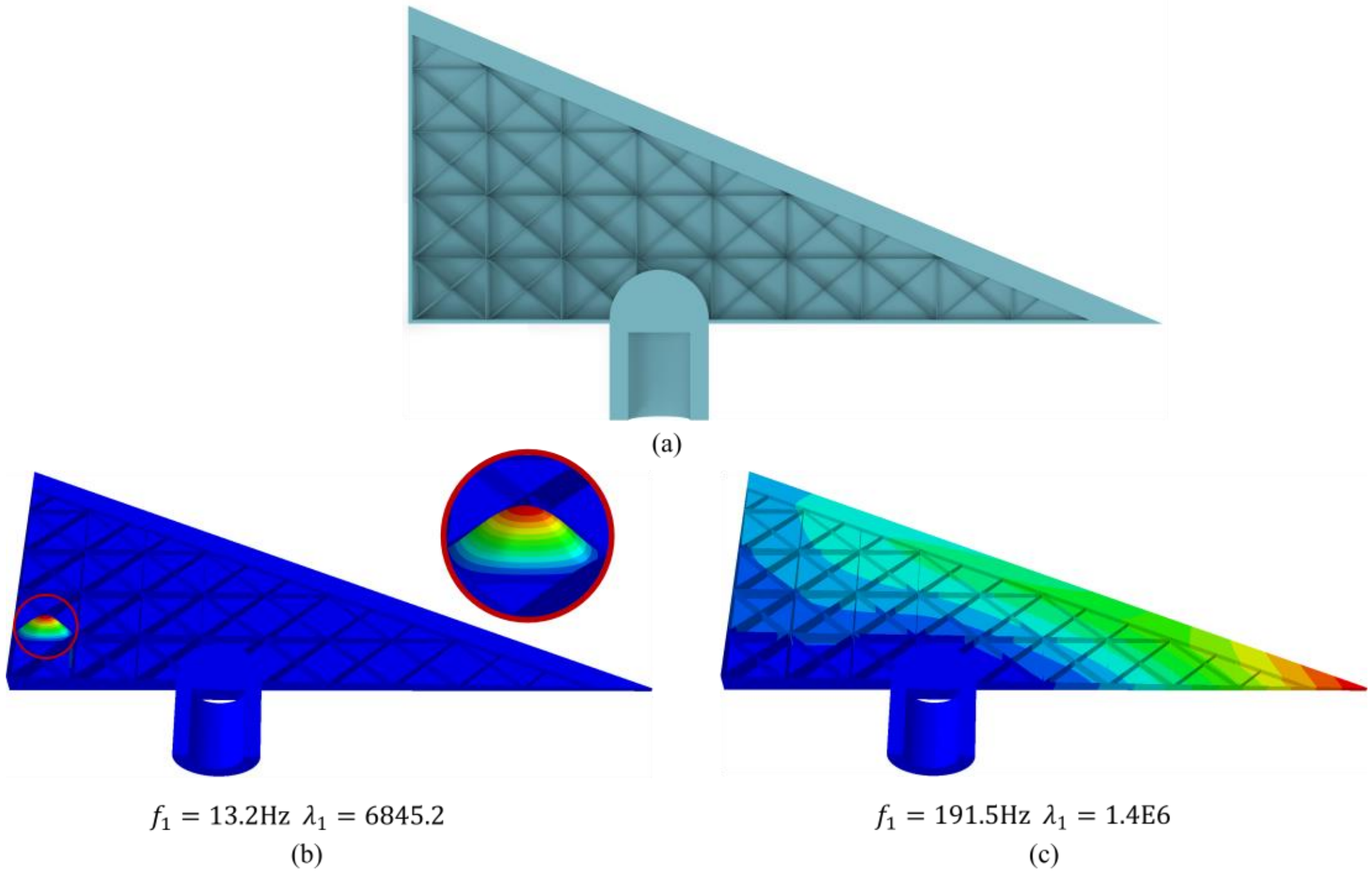

Fig. 11. Validation of density penalty scheme: (a) initial layout (displaying thickness); (b) spurious mode without using density penalty scheme; (c) normal vibration mode obtained using density penalty scheme.

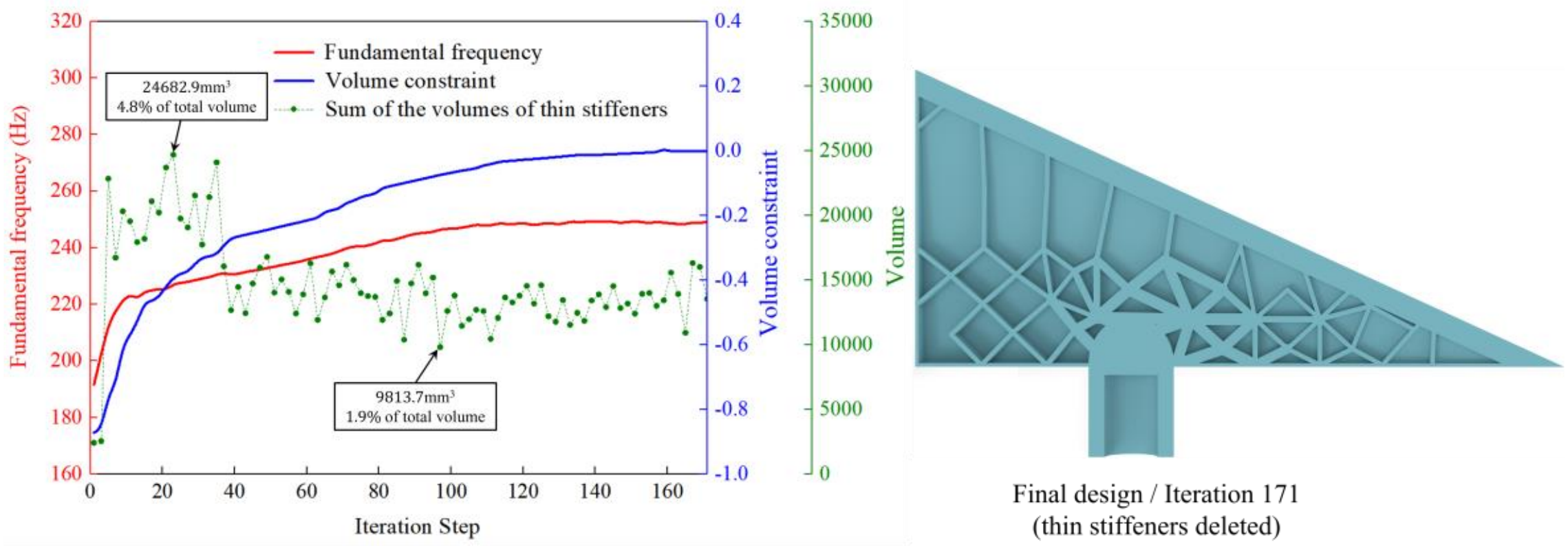

(a) (b)

Fig. 12. (a) Iteration history of objective function and constraint values; (b) final design of the triangular air rudder.

-unterpart removing the stiffeners with thickness smaller than 4mm , those vibration modes are consistently out-of-plane bending pattern with no local vibrations of thin stiffeners. Furthermore, the relative difference of the fundamental frequency (with and without removing stiffeners) is within 0.3%. This suggests that the removal of thin stiffeners does not significantly affect the fundamental frequency and vibration modes of the triangular rudder during the optimization process, validating the effectiveness of density penalization scheme on suppressing spurious modes.

Table 4 Comparison of the different designs for rudder structure.

| | Optimized design | Original design |
|---|---|---|
| Stiffener distribution (displaying thickness) | 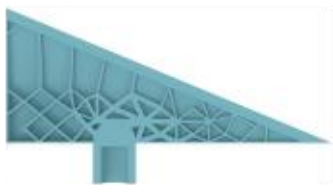 | 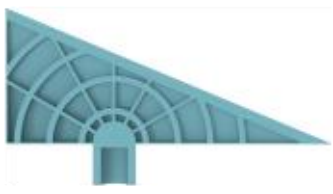 |
| Vibration mode | 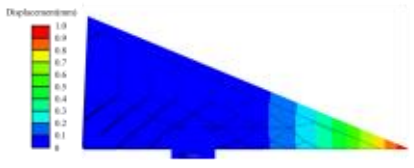 | 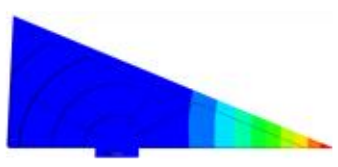 |
| Fundamental frequency (Hz) | 249.2 | 217.0 |
| Total mass (kg) | 5.28 | 5.28 |

Table 5 Fundamental frequency and vibration modes of some intermediate configurations.

| | fundamental frequency (with thin stiffeners) | fundamental frequency (thin stiffeners deleted) |
|---|---|---|
| Iteration 11 | 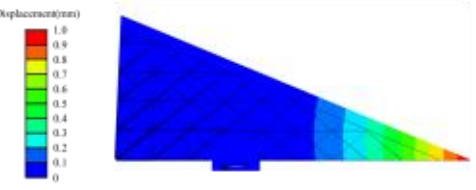 $f_1 = 223.8$Hz | 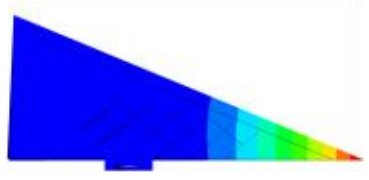 $f_1 = 223.2$Hz |
| Iteration 80 | 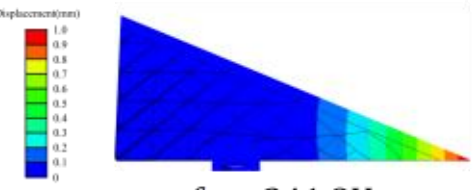 $f_1 = 241.8$Hz | 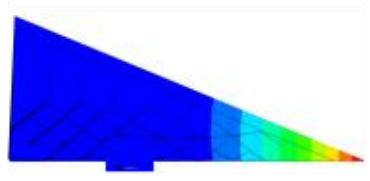 $f_1 = 242.3$Hz |
| Iteration 171 | 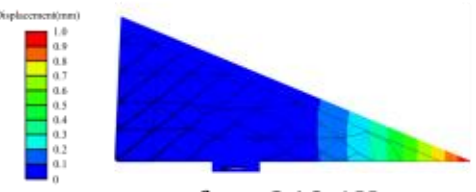 $f_1 = 249.1$Hz | 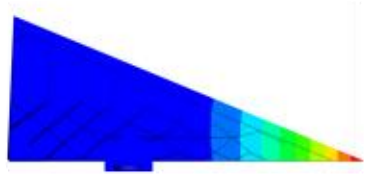 $f_1 = 249.2$Hz |

## 5. Conclusions

This paper presents an efficient design method for air rudders commonly used for controlling aircraft attitude, and the following main achievements are obtained:

(1) Based on the MMC optimization design framework, the stiffeners within the irregular enclosed domain are accurately described and adaptively modeled using explicit parameters and shell elements.

(2) For design objectives of minimizing structural compliance and maximizing fundamental frequency, analytical sensitivities are derived and validated by numerical examples. To suppress the spurious modes and satisfy manufacturing requirement, penalization schemes about the thickness and density of stiffeners are proposed.

(3) The explicit description of the stiffener in the optimized air rudder enables CAD modelling, thereby eliminating the tedious process of manual identification and further optimizing of stiffeners in implicit topology optimization approach.

The present work can be extended for the intelligent design of other aerospace structures beyond air rudders. Ongoing research also focuses on the buckling optimization and thermomechanical design optimization of thin-walled structures in aerospace engineering.

## Acknowledgements

This research was funded by the National Natural Science Foundation, China (12002073, 12372122), the National Key Research and Development Plan, China (2023YFB3309104), the Science Technology Plan of Liaoning Province (2023JH2/101600044), the Liaoning Revitalization Talents Program (XLYC2001003) and 111 Project, China (B14013).

## References

1. Bendsoe MP, Kikuchi N. Generating optimal topologies in structural design using a homogenization method. *Comput Methods Appl Mech Eng* 1988;**71**:197–224.
2. Yıldız AR, Kılıçarpa UA, Demirci E, et al. Topography and topology optimization of diesel engine components for lightweight design in the automotive industry. *Mater Test* 2019;**61**(1):27–34.
3. Boccini E, Furferi R, Governi L, et al. Toward the integration of lattice structure-based topology optimization and additive manufacturing for the design of turbomachinery components. *Adv Mech Eng* 2019;**11**(8):1–14.
4. Dbouk T. A review about the engineering design of optimal heat transfer systems using topology optimization. *Appl Therm Eng* 2017;**112**:841–854.
5. Remouchamps A, Bruyneel M, Fleury C, et al. Application of a bilevel scheme including topology optimization to the design of an aircraft pylon. *Struct Multidiscip Optim* 2011;**44**(6):739–750.
6. Xie GN, Qi W, Zhang WH, et al. Optimization design and analysis of multilayer lightweight thermal protection structures under aerodynamic heating conditions. *J Therm Sci Eng Appl* 2013;**5**(1):011011.
7. Hou J, Zhu JH, Fei H, et al. Stiffeners layout design of thin-walled structures with constraints on multi-fastener joint loads. *Chin J Aeronaut* 2017; **30**(4): 1441–1450.
8. Li SY, Wei HK, Yuan SQ, et al. Collaborative optimization design of process parameter and structural topology for laser additive manufacturing. *Chin J Aeronaut* 2023; **36**(1):456–467.
9. Zhao X, Zhang WH, Ying Z, et al. Multiscale topology optimization using feature-driven method. *Chin J Aeronaut* 2020; **33**(2):621–633.
10. Song LL, Gao T, Tang L, et al. An all-movable rudder designed by thermo-elastic topology optimization and manufactured by additive manufacturing. *Comput Struct* 2021;**243**:106405.
11. Wang C, Zhu JH, Wu MQ, et al. Multi-scale design and optimization for solid-lattice hybrid structures and their application to aerospace vehicle components. *Chin J Aeronaut* 2021; **34**(5):386–398.
12. Zhu JH, Zhao YB, Zhang WH, et al. Bio-inspired feature-driven topology optimization for rudder structure design. *Eng Sci* 2018;**5**(2):46–55.
13. Gu XJ, Yang KK, Wu MQ, et al. Integrated optimization design of smart morphing wing for accurate shape control. *Chin J Aeronaut* 2021;**34**(1): 135–147.
14. Bendsøe MP. Optimal shape design as a material distribution problem. *Struct Optim* 1989;**1**(4):193–202.
15. Mlejnek HP. Some aspects of the genesis of structures. *Struct Optim* 1992;**5**(1–2):64–69.
16. Bendsoe MP, Guedes JM, Haber RB, et al. An analytical model to predict optimal material properties in the context of optimal structural design. *J Appl Mech* 1994;**61**(4):930–937.
17. Wang MY, Wang XM, Guo DM. A level set method for structural topology optimization. *Comput Methods Appl Mech Eng* 2003;**192**(1):227–246.
18. Allaire G, Jouve F, Toader AM. Structural optimization using sensitivity analysis and a level-set method. *J Comput Phys* 2004;**194**(1):363–393.
19. Xie YM, Steven GP. A simple evolutionary procedure for structural optimization. *Comput Struct* 1993;**49**(5):885–896.
20. Xie YM, Steven GP. Optimal design of multiple load case structures using an evolutionary procedure. *Eng Comput* 1994;**11**(4):295–302.
21. Guo X, Cheng GD. Recent development in structural design and optimization. *Acta Mech Sin* 2010;**26**(6):807–823.
22. Sigmund O, Maute K. Topology optimization approaches. *Struct Multidiscip Optim* 2013;**48**(6):1031–1055.
23. Wang C, Zhao Z, Zhou M, et al. A comprehensive review of educational articles on structural and multidisciplinary optimization. *Struct Multidiscip Optim* 2021;**64**(5):2827–2880.
24. Du ZL, Jia YB, Chung HY, et al. Analysis and optimization of thermoelastic structures with tension–compression asymmetry. *Int J Solids Struct* 2022;**254–255**:111897.
25. Zhu JH, Zhang WH, Xia L. Topology optimization in aircraft and aerospace structures design. *Arch Comput Methods Eng* 2016;**23**(4):595–622.
26. Shi GH, Jia YB, Hao WY, et al. Optimal design of rudder structures based on data-driven method. *Chin J Theor Appl Mech* 2023;**55**(11):2577–2587.
27. Chen JY, Xin KH, Gu XL, et al. Topology optimization and design for additive manufactured supporting structure of vehicle rudder. *J Phys Conf Ser* 2021;**2065**(1):012022.
28. Feng ZJ, Du JB, Zhang HX, et al. Structural design flow of typical aircraft components based on topology optimization. *IOP Conf Ser Mater Sci Eng* 2020;**892**(1):012029.
29. Guo X, Zhang WS, Zhong WL. Doing topology optimization explicitly and geometrically—a new moving morphable components based framework. *J Appl Mech* 2014;**81**:081009.
30. Zhang WS, Yuan J, Zhang J, et al. A new topology optimization approach based on Moving Morphable Components (MMC) and the ersatz material model. *Struct Multidiscip Optim* 2016;**53**(6):1243–1260.
31. Du ZL, Cui TC, Liu C, et al. An efficient and easy-to-extend Matlab code of the Moving Morphable Component (MMC) method for three-dimensional topology optimization. *Struct Multidiscip Optim* 2022;**65**(5):158.
32. Jiang XD , Liu C, Zhang SH, et al. Explicit topology optimization design of stiffened plate structures based on the Moving Morphable Component (MMC) method. *Comput Model Eng Sci* 2023;**135**(2):809–838.

33. Jiang XD, Liu C, Du ZL, et al. A unified framework for explicit layout/topology optimization of thin-walled structures based on Moving Morphable Components (MMC) method and adaptive ground structure approach. *Comput Methods Appl Mech Eng* 2022;**396**:115047.
34. Jiang XD, Huo WD, Liu C, et al. Explicit layout optimization of complex rib-reinforced thin-walled structures via computational conformal mapping (CCM). *Comput Methods Appl Mech Eng* 2023;**404**:115745.
35. Huo WD, Liu C, Liu YP, et al. A novel explicit design method for complex thin-walled structures based on embedded solid moving morphable components. *Comput Methods Appl Mech Eng* 2023;**417**:116431.
36. Guo X, Zhou JH, Zhang WS, et al. Self-supporting structure design in additive manufacturing through explicit topology optimization. *Comput Methods Appl Mech Eng* 2017;**323**:27–63.
37. Komkov V, Choi KK, Haug EJ. *Design Sensitivity Analysis of Structural Systems.* Academic Press; 1986.
38. Laporte E, Tallec PL. *Numerical Methods in Sensitivity Analysis and Shape Optimization.* Springer Science & Business Media;2002.
39. Seyranian AP, Lund E, Olhoff N. Multiple eigenvalues in structural optimization problems. *Struct Optim* 1994;**8**(4):207–227.
40. Bendsoe MP, Sigmund O. *Topology Optimization: Theory, Methods, and Applications.* Springer Science & Business Media;2003.
41. Du JB, Olhoff N. Topological design of freely vibrating continuum structures for maximum values of simple and multiple eigenfrequencies and frequency gaps. *Struct Multidiscip Optim* 2007;**34**(2):91–110.
42. Liu Y, Shimoda M. Non-parametric shape optimization method for natural vibration design of stiffened shells. *Comput Struct* 2015;**146**:20–31.
43. Svanberg K. The method of moving asymptotes—a new method for structural optimization. *Int J Numer Methods Eng* 1987;**24**(2):359–373.